\documentclass[12pt,centertags,oneside]{amsart}

\usepackage{amsmath,amstext,amsthm,amssymb,amsfonts,amscd}
\usepackage{mathrsfs}
\usepackage[colorlinks=true]{hyperref}
\usepackage{graphicx}
\usepackage{todonotes}
\usepackage{tocvsec2}
\usepackage{typearea}
\usepackage[shortlabels]{enumitem}
\usepackage{tikz-cd}

\usepackage{fouriernc} 

\usepackage[a4paper,width=14.66cm,top=2.52cm,bottom=2.52cm]{geometry}

\usepackage[capitalise]{cleveref}
 
\crefname{thm}{}{}
\crefname{prop}{}{}
\crefname{re}{}{}
\crefname{def}{}{}
\crefname{equation}{}{}

\theoremstyle{plain}
\newtheorem{lem}{Lemme}[section]
\newtheorem{thm}[lem]{Theorem}
\newtheorem{prop}[lem]{Proposition}
\newtheorem{cor}[lem]{Corollary}

\theoremstyle{definition}
\newtheorem{defin}[lem]{Definition}
\newtheorem{exa}[lem]{Example}

\theoremstyle{remark}
\newtheorem{re}[lem]{Remark}

\numberwithin{equation}{section}
\numberwithin{figure}{section}

\newcommand{\cO}{\mathcal{O}}

\newcommand{\bC}{\mathbf{C}}

\DeclareMathOperator{\Sp}{\mathrm{Sp}}

\renewcommand{\Re}{\mathrm{Re}\,}
\DeclareMathOperator{\im}{\mathrm Im}
\DeclareMathOperator{\End}{\mathrm End}

\begin{document}

\title[Geometric Zabrodin-Wiegmann conjecture]{Geometric Zabrodin-Wiegmann conjecture for integer Quantum Hall states}

\author{Shu Shen}
\address{Institut de Math\'ematiques de Jussieu-Paris Rive Gauche,
	Sorbonne Universit\'e, Case Courrier 247, 4 place Jussieu, 75252 Paris Cedex 5, France.}
\email{shu.shen@imj-prg.fr}

\author{Jianqing Yu}
\address{School of Mathematical Sciences, University of Science and Technology of China, 96 Jinzhai Road, Hefei, Anhui 230026, P. R. China.}
\email{jianqing@ustc.edu.cn}

\thanks {}

\makeatletter
\@namedef{subjclassname@2020}{%
	\textup{2020} Mathematics Subject Classification}
\makeatother 

\begin{abstract} 
The purpose of this article is to show a geometric version of Zabrodin-Wiegmann conjecture for an integer quantum Hall state.
Given an effective reduced divisor on a compact connected Riemann surface, using the canonical holomorphic section of the associated canonical line bundle as well as certain initial data and local normalisation data,
we construct a canonical non-zero element in the determinant line of the cohomology of the $ p$-tensor power of the line bundle.
 
When endowed with proper metric data, the square of the $ L^{2} $-norm of our canonical element is the partition function associated to an integer quantum Hall state.
We establish an asymptotic expansion for the logarithm of the partition function when $ p\to +\infty$.
The constant term of this expansion includes the holomorphic analytic torsion and matches 
a geometric version of Zabrodin-Wiegmann's prediction.

Our proof relies on Bismut-Lebeau's embedding formula for the Quillen metrics, Bismut-Vasserot and Finski's asymptotic expansion for the analytic torsion associated to the higher tensor product of a positive Hermitian holomorphic line bundle.
\end{abstract}


\date{\today}

\dedicatory{}

\maketitle
\tableofcontents
\settocdepth{subsection}
\numberwithin{equation}{section}

\section*{Introduction}
This paper studies the mathematical models of the quantum Hall states and the Quillen metric.

The Quillen metric was introduced by Quillen \cite{Quillendeter} for the determinant line of the cohomology of a holomorphic vector bundle on a compact Riemann surface. 
Its key properties were developed in the works of Bismut-Gillet-Soulé \cite{BGS1,BGS2,BGS3} and Bismut-Lebeau \cite{BL91}.

In the 1990s, the Quillen metric found applications in modeling the quantum Hall effect, helping to understand the quantization of transport coefficients \cite{AvronZograf94,AvronSeilerZograf95,Levay95,Levay97}. 
More recently, in \cite{KlevtsovMaMarinescuWiegmann17}, 
Klevtsov, Ma, Marinescu, and Wiegmann  used the anomaly formula of Bismut-Gillet-Soulé \cite{BGS3} to study the integer quantum Hall effect.

Building on these developments, we further investigate Bismut-Lebeau's embedding formula \cite{BL91} and rigorously confirm a geometric version of Zabrodin-Wiegmann conjecture \cite{ZabrodinWiegmann06} for the partition function associated to an integer quantum Hall state.

\subsection{Physical background}
The quantum Hall effect provides a paradigmatic example of collective quantum phenomena in two-dimensional electron systems subjected to a strong perpendicular magnetic field.

In the integer quantum Hall regime, the many-body ground state can be described explicitly in terms of single-particle eigenfunctions of the lowest Landau level.
More precisely, for a uniform magnetic field, parameterized by a positive constant $ B$, the quantum dynamics of an electron in the plane $ \mathbf{R}^{2} $ is governed by the magnetic Laplacian, 
\begin{align}\label{eq:hsy2}
 H_{B} = -\left(\frac{\partial }{\partial x} + i \frac{B}{2}  y\right)^{2} - \left(\frac{\partial }{\partial y} -i \frac{B}{2}   x\right)^{2}. 
\end{align}
It is well-known that $ H_{B} $ is essentially self-adjoint and has spectrum (Landau levels) given by $ \left(2\mathbf{N} +1\right)B$.
In complex coordinates $ z= x+i y$, the eigenspace corresponding to the lowest Landau level admits a basis given by
\begin{align}\label{eq:nxnd}
 s_{j} \left(z\right)= z^{j} e^{-\frac{B}{4}  \left| z \right|^{2} } \text{ with } j\in \mathbf{N}.
\end{align}
The radial density associated to $ s_{j} $ is localized around the radius\footnote{The measure associated to $ s_{j} $ is given by $ \left| s_{j}\left(z\right) \right|^{2} dxdy = r^{2j+1} e^{-Br^{2} /2} drd\theta$. The radial density $ r^{2j+1} e^{-Br^{2} /2}$ takes its maximal value at $ r_{j} = \sqrt{\frac{2j+1}{B} }$.} $ r_{j} = \sqrt{\frac{2j+1}{B} }$, with a fast decay away from this peak.

For a finite droplet of $ N$  electrons at filling factor $ \nu = 1$, the occupied orbitals are those with  $j=  0, 1,2, \ldots, N-1$ forming a compact ``quantum Hall droplet'' of radius $ R\propto \sqrt{N/B}$.
In typical scaling limits one takes $ R\propto 1$, which corresponds to $ N \propto B$.
To satisfy the Pauli exclusion principle, the many-body wave function of $ N$ electrons is antisymmetric, and is then the Slater determinant of these lowest orbitals.
This wave function simplifies to the Vandermonde form, 
\begin{align}\label{eq:irvd}
  \Delta_{N}  \left(z_{1},\ldots, z_{N} \right) \prod_{j= 1}^{N} \exp\left(-\frac{B}{4} \left| z_{j}  \right|^{2}  \right),
\end{align}
where $ \Delta_{N}  \left(z_{1},\ldots, z_{N} \right) = \prod_{1\le i<j \le N} \left(z_{j} -z_{i} \right)  $ is the Vandermonde determinant.

In the fractional quantum Hall regime, electron–electron interactions play a decisive role.
In this setting, the eigenfunctions of the corresponding magnetic $ N$-body Hamiltonian are not explicitly known. 
To describe the fractional quantum Hall effect, Laughlin \cite{Laughlin81} introduced the trial states
\begin{align}\label{eq:zoku}
 \left( \Delta_{N}  \left(z_{1},\ldots, z_{N} \right) \right)^{\beta } \prod_{j= 1}^{N} \exp\left(-\frac{B}{4} \left| z_{j}  \right|^{2}  \right),
\end{align}
where \(\beta\) is an odd positive integer.

If $ W : \mathbf{C} \to \mathbf{R} $ is some suitable function on $ \mathbf{C} $, put 
\begin{align}\label{eq:ciyh}
 Z_{\beta, N} = \frac{1}{N!} \int_{\mathbf{C}^{N} }  \left|  \Delta_{N}  \left(z_{1},\ldots, z_{N} \right) \right|^{2\beta }\prod_{j= 1}^{N} \exp\left(-N W\left(z_{j} \right) \right)dx_{j}dy_{j}.
\end{align}
The asymptotic analysis as $ B\propto N\to +\infty$  for the logarithm of the norms of \eqref{eq:irvd}, \eqref{eq:zoku}, or their natural normalised generalisation  $ \log Z_{\beta, N} $ plays a crucial role in the understanding of both the integer ($ \beta=1$) and fractional ($\beta\neq1$) quantum Hall regimes.

The formal expression \eqref{eq:ciyh} coincides with the partition function of a statistical ensemble of $ N$ two-dimensional Coulomb charges (2D Dyson gas) at temperature $T= 1/\beta$ in the external potential $ W$.
For a fixed $ \beta $, Zabrodin-Wiegmann \cite{ZabrodinWiegmann06} studied the asymptotic expansion of $ \log Z_{\beta, N} $ as $ N \to + \infty $.
Using Ward identity, they predicted the expansion \cite[(1.2) and Footnote 6]{ZabrodinWiegmann06},
\begin{align}\label{eq:pgtp}
 \log Z_{\beta, N}  = N^{2} F_{0 } + \frac{\beta -2}{2} N\log N + N  F_{1/2} + F_{1 } + \mathcal{O} \left(1/N\right),
\end{align}
where the coefficients $ F_{0}, F_{1/2}, F_{1} $ are  explicitly given in \cite[(3.5)-(3.9), (4.6)-(4.9)]{ZabrodinWiegmann06}.

In a two-dimensional quantum Hall system, as $ N \to +\infty$, the $ N$ electrons condense into a droplet occupying a domain $ D \subset \mathbf{C} $ determined by the function $ W$.
The leading term $ F_0$ encodes the bulk electrostatic energy.
The subleading corrections $F_{1/2} $ and $ F_1$ reflect 
the geometry of the droplet.
A key feature of Zabrodin and Wiegmann's prediction is that $ F_{1} $, in addition to classical terms, includes global spectral invariants such as the relative regularized determinant of the Laplacian on $\mathbf{C} \backslash D $.

To further explore the geometric and topological aspects of quantum Hall states, we place the system on a torus or, more generally, on a compact Riemann surface endowed with a holomorphic vector bundle \cite{AvronRavehZur88,Klevtsov14,KlevtsovMaMarinescuWiegmann17,KlevtsovZvonkine22,BurbanKlevtsov24}.
This setup incorporates a curved background metric, spatially varying magnetic fields, spin line bundles encoding the geometric spin structure, and flat line bundles associated with nontrivial $ 1$-cycles that support Aharonov-Bohm flux insertions (see Footnotes \ref{f2} and \ref{f3}).
The partition function in this regime reveals new physical phenomena, including the Wen-Zee term \cite{WenZee92} \cite[(3.21)]{Klevtsov16}.
The term $ F_{1} $ in the corresponding large $ N$ expansion captures quantum and gravitational anomalies that are sensitive to the global geometry and topology of the surface.
Moreover, when the complex structure of the Riemann surface or the holomorphic vector bundle varies, or when the metric changes, the resulting variations of the partition function are related to the so-called adiabatic transport coefficients \cite{Avron95}.
Altogether, these aspects connect to index theory and complex geometry, offering deeper insight into the interplay between physical phenomena and rich mathematical structures.

\subsection{Mathematical background}
Let us formulate a geometric version for the Zabrodin-Wiegmann conjecture on a Riemann surface within a precise mathematical setting.
We use standard terminology on Riemann surfaces and refer the reader to Section \ref{S:Prel} for precise definitions.

Let $ X$ be a compact connected Riemann surface of genus $ g$.
Let $ L$ be a positive line bundle\footnote{\label{f2}In physical applications, $ L$ is typically taken to have degree one.
Given a Hermitian metric on $ L$, the curvature of the Chern connection on $ L^{p} $ is a magnetic field with total flux $ p$ through the Riemann surface $ X$, generated by a distribution of magnetic monopole charges.} on $ X$.
Let $ E$ be a holomorphic vector bundle\footnote{\label{f3}One often takes $ E=  K_{X}^{s} \otimes F$, where $ s\in \frac{1}{2} \mathbf{Z} $ is the spin, $ K_{X} $ is the canonical line bundle, and $ F$ is a unitary flat line bundle.
The latter amounts to choosing a unitary representation of the fundamental group $ \pi_{1} \left(X\right)$, whose values give the holonomies around non-trivial $ 1$-cycles of $ X$.
These holonomies encode the Aharonov–Bohm fluxes, with each cycle carrying a phase determined by the representation.
} on $ X$.
For $ p\in \mathbf{N} $,  let $ H \left(X,L^{p}  \otimes E\right)= H^{0 }  \left(X,L^{p}  \otimes E\right)\oplus H^{1} \left(X,L^{p}  \otimes E\right)$ be the cohomology of the sheave of holomorphic sections of $ L^{p} \otimes E$.
Let 
\begin{align}\label{eq:hg33}
 \lambda_{p} \left(E\right)= \det H^{0 } \left(X,L^{p}  \otimes E\right)\otimes \left( \det H^{1 } \left(X,L^{p}  \otimes E\right)\right)^{-1} 
\end{align}
be the determinant line of the cohomology of $ L^{p} \otimes E$.


Assume $ p$ is sufficiently large such that $ \deg \left(L^{p}\otimes E\right) > \deg K_{X} $.
Then,  $ H^{1} \left(X,L^{p}  \otimes E\right)$ vanishes and $ \dim H^{0 } \left(X,L^{p}  \otimes E\right) $ is given the Riemann-Roch formula, 
\begin{align}\label{eq:fr2d}
 N_{p} = \left(1-g +p \deg\left(L\right) \right){\rm rk}\left(E\right) + {\rm deg}\left(E\right).
\end{align}
In this case, 
\begin{align}\label{eq:trgm}
 \lambda_{p} \left(E\right)= \det H^{0 } \left(X,L^{p}  \otimes E \right).
\end{align}

Let $ \left(s_{p}^{i}  \right)_{1\le i \le N_{p} } $ be a basis of $ H^{0 } \left(X,L^{p}  \otimes E\right)$.
The Slater determinant $ \det \left(s^{i}_ {p}  \left(z_{j} \right) \right)$ of this family is a section of $ \left(L^{p} \otimes E \right)^{ \boxtimes N_{p} } $ over $ X^{N_{p} } $  defined using the standard determinant formula for matrices  (see \eqref{eq:wmbo}).

Let  $ \omega^{X} $ be a K\"ahler form on $ X$, let $ h^{L}, h^{E} $ be Hermitian metrics  on $ L$, $ E$.
Assume $ \left(\omega^{X},h^{L}  \right)$ satisfies the prequantization condition (see \eqref{eq:yagj}). 

Let $ \left| \det \left(s^{i}_{p}  \left(z_{j} \right) \right) \right|_{h^{L}, h^{E} }  $ be the norm of $ \det \left(s^ {i}_ {p}  \left(z_{j} \right) \right)$ with respect to the  Hermitian metric on $ \left(L^{p} \otimes E \right)^{ \boxtimes N_{p} }$ induced by $ h^{L}, h^{E} $.
Let $ dv_{X^{N_{p}  } } $ be the normalised (see \eqref{eq:epzu}) volume form on $ X^{N_{p} } $.
For $ \beta > 0  $, the partition function is now defined by 
\begin{align}\label{eq:erow}
 Z_{p } = \frac{1}{N_{p} !} \int_{X^{N_{p} }}^{} \left| \det \left(s^{i}_ {p}  \left(z_{j} \right) \right) \right|^{2 \beta }_{h^{L}, h^{E} }  dv_{X^{N_{p} } }.
\end{align}
When $ \beta = 1$, $ Z_{p }$ is the square of the $ L^{2} $-norm of the integer quantum Hall state $ \frac{1}{\sqrt{N_{p} !}} \det \left(s^{i}_ {p}  \left(z_{j} \right) \right)$.

A natural question is whether the same phenomenon observed in \eqref{eq:pgtp} occurs for $ \log Z_{p } $ as $ p\to +\infty$.
Additionally, one may ask whether global spectral invariants, such as the zeta-regularized determinant of certain Laplacians, appear in the constant term of the expansion.
Following a suggestion by Klevtsov, we call this the geometric Zabrodin-Wiegmann conjecture.

As in \cite{Klevtsov14,KlevtsovMaMarinescuWiegmann17}, we now restrict ourself to the case $ \beta = 1$.\footnote{See \cite{Klevtsov16} for more details on this conjecture for general $ \beta\in \mathbf{R}_{+}^{*} $.}
In this case, it has been noted in the references above that the partition function $ Z_{p }$ is related to the $ L^{2} $-norm on the determinant line $\lambda_{p} \left(E\right)$.
Indeed, put 
\begin{align}\label{eq:ngav}
 s_{p} = s_{p}^{1}  \wedge s_{p}^{2}\wedge  \cdots \wedge s^ {N_{p}}_ {p}.
\end{align}
Then, $ s_{p} $ is a non-zero element in $ \lambda_{p} \left(E\right)$.
If we equip $ H^{0 } \left(X,L^{p} \otimes E\right)$ the $ L^{2} $-metric induced from Hodge theory, then $ \lambda_{p} \left(E\right)$ inherits an induced norm $ \left| \cdot \right|_{\lambda_{p} \left(E\right)} $.
We have the classical formula (\cite[(2.6)]{Klevtsov14}, see also Proposition \ref{prop:apm1}), 
\begin{align}\label{eq:jxwb}
	Z_{p } = \left| s_{p}  \right|^{2}_{\lambda_{p} \left(E\right)}.  
\end{align}
 
The behavior of $ Z_{p }  $ depends on the choice of $ s_{p}  \in \lambda_{p} \left(E\right)$.
In our paper, using some extra information from an effective reduced divisor $ D$, we construct canonically a non-zero element $ s_{p}\in \lambda_{p} \left(E\right)$, and confirm the geometric Zabrodin-Wiegmann conjecture for our section $ s_{p} $.
We also show that our $ s_{p} $ is compatible with the previously known constructions in the cases $ g= 0 $ and $ 1$.

\subsection{Main results}
Let us state our main results.
Let $ D \subset X$ be an effective reduced divisor of $ X$.
Let $ L$ be the associated positive holomorphic line bundle, and let $ s_{D} $ be the corresponding holomorphic section of $ L$.

Put 
\begin{align}\label{eq:mo1w}
& \lambda \left(L_{|D} \right)= \bigotimes_{z\in D}  L_{z},&\lambda \left(E_{|D} \right)= \bigotimes_{z\in D} \det \left(E_{z}\right).
\end{align}
To be precise, the first is a tensor product of even lines, whereas the second is a tensor product of signed lines (see Section \ref{sec:1c1y}).

\begin{thm}\label{thm:e2jm}
 Given an initial data $ s_{0 }\in \lambda_{0}\left(E\right) $ and some non-zero local normalisation data $ s^{L}_{D}\in \lambda \left(L_{|D} \right) $,  $ s^{E}_{D}\in \lambda \left(E_{|D}\right)  $, one can use  $ s_{D} $ to canonically associate,  for each $ p\in \mathbf{N} $,  a non-zero element $ s_{p} $ in $ \lambda_{p}\left(E\right) $.
 
Moreover, with metric data $\left( \omega^{X}, h^{L}, h^{E} \right)$ such that $\left(\omega^{X}, h^{L} \right) $ satisfies the prequantization conditions, the associated partition function $ Z_{p} \left(s_{0 },s^{L}_{D} ,  s^{E} _{D} \right)$ has the asymptotic expansion as $ p\to +\infty$, 
 \begin{align}\label{eq:zjj1}
\log	Z_{p} \left(s_{0 },s^{L}_{D} ,  s^{E} _{D} \right) = a_{2} p^{2}+b_{1} p\log p +a_{1}p+b_{0 } \log p+ a_{0 }  +\mathcal{O} \left(\frac{\log\left(p\right)}{p}  \right).
 \end{align} 
 All the constants above are explicitly determined in \eqref{eq:koff} and \eqref{eq:dbkb}.
 In particular, 
 \begin{align}\label{eq:aaaz}
	a_{0 } = \log \left|s_{0 } \right|_{\lambda_{0 }  \left(E\right)} ^{2}	+2\tau\left(\omega^{X}, h^{E} \right)  {  -}  \left(\zeta'\left(-1\right)+\frac{\log\left(2\pi \right)}{12}+\frac{7}{24}\right){\rm rk}\left(E\right)\chi\left(X\right)
	{  -} \frac{1}{2}\deg\left(E\right),
 \end{align}
 where $ \tau\left(\omega^{X}, h^{E} \right)$ is the Ray-Singer holomorphic torsion of $ E$,  $ \chi\left(X\right)= 2-2g$ is the Euler characteristic of $ X$, and $ \zeta $ is the Riemann zeta function.
\end{thm}

The term $ a_{0 } $ is independent of $ s_{D} $ and of the local normalisation data $ \left(s^{L}_{D}, s^{E}_{D}  \right)$.
The holomorphic torsion of $ \tau \left(\omega^{X}, h^{E} \right)$ is precisely the global spectral invariant as predicted by Zabrodin-Wiegmann \cite{ZabrodinWiegmann06}.
Recall that if $ \Box^{E} $ is the Kodaira Laplacian of $ E$ (see \eqref{eq:ztwh}), we have  
\begin{align}\label{eq:zcxu}
	\tau \left(\omega^{X}, h^{E} \right) =-\frac{1}{2}  \log {\rm det}_{\zeta } \left(\Box^{E} _{| \Omega^{0,1}\left(X,E\right)} \right)= - \frac{1}{2}  \log {\rm det}_{\zeta } \left(\Box^{E} _{| C^{\infty} \left(X,E\right)} \right),
\end{align}
where $ \det_{\zeta }  $ is the zeta-regularized determinant.

\begin{re}\label{re:jzrq}
	In Theorem \ref{thm:wlcn}, we obtain a more general version of Theorem \ref{thm:e2jm}.
	Indeed, we will establish a full asymptotic expansion for $ \log Z_{p}\left(s_{0 },s^{L}_{D} ,  s^{E} _{D} \right) $ when $ p\to +\infty$ for Hermitian metric $ h^{L} $ with positive curvature, without requiring the prequantization condition.
	Furthermore, the constants $ a_{1 }, a_{2 }, b_{0 }, b_{1} $ are explicitly evaluated in this general case. Now, $ a_{1} $ depends also on the ratio of the curvature of $ h^{L} $ to the K\"ahler form  $ \omega^{X} $, which shows how far it is from being prequantized.
	However, to evaluate the constant $ a_{0 } $, we need to impose the prequantization condition.
	Otherwise, the derivations of the default term mentioned above would also appear.	 
\end{re}

\subsection{Constructions of $ s_{p}\in \lambda_{p} $ }\label{isConstr}
The canonical holomorphic section $ s_{D}$ of $ L$ is the key in our construction of $ s_{p} $. 
It plays the role of the link between $ s_{p-1} $ and $ s_{p} $.

More precisely, if $ \iota : D \to X$ is the natural embedding, if $ \mathcal{O}_{X} \left(L^{p} \otimes E\right)$ is the sheave of holomorphic sections of $ L^{p} \otimes E$, and 
if $ \left(L^{p} \otimes E\right)_{|D} $ is the obvious sheave defined on $ D$, we have the exact sequence of sheaves on $ X$, 
\begin{equation}\label{dia:pgtl0 }
	\begin{tikzcd}
		0 & \mathcal{O}_{X} \left(L^{p-1} \otimes E\right) &	\mathcal{O}_{X} \left(L^{p} \otimes E\right) & \iota_{*}\left(\left(L^{p} \otimes E\right)_{|D}\right) & 0,
		 \arrow["", from=1-1, to=1-2]
		 \arrow["s_{D} ", from=1-2, to=1-3]
		 \arrow["r_{D} ", from=1-3, to=1-4]
		 \arrow["", from=1-4, to=1-5]
	\end{tikzcd} 
\end{equation}
where the first map is the multiplication by $ s_{D}$, and the second map is the restriction.

The theory of determinant of Grothendieck and Knudsen-Mumford \cite{KM76} provides us a canonical isomorphism,  
\begin{align}\label{eq:3ko2}
 \sigma^{p-1}_{p} : \lambda_{p-1} \left(E\right) \otimes \det \left(\left(L^{p} \otimes E\right)_{|D}\right) \to \lambda_{p} \left(E\right).
\end{align}
Using $ \sigma_{p}^{p-1} $ and a non-zero element in $ \det \left(\left(L^{p} \otimes E\right)_{|D}\right)$ determined by the local normalisation data $\left( s^{L}_{D}, s^{E}_{D}\right) $, we can define inductively a non-zero element $ s_{p}\in  \lambda_{p} \left(E\right)$.

We note the similarity of this construction with Morse theory.
If one has a Morse function on a real smooth manifold, the topology of the manifold can be recovered by shifting the level sets of the Morse  function and crossing the critical points one by one. 
Each step involves only a simple local model.

This analogy is not merely philosophical. In fact, the Milnor metric, for a Morse-Smale flow is constructed in this manner \cite[Section 2]{ShenYuMorseSmale}. Our construction of $ s_{p} $  is inspired by this approach.

\begin{re}\label{re:xbot}
To construct the canonical isomorphism $ \sigma_{p}^{p-1} $ in \eqref{eq:3ko2}, we must use the full formalism developed in \cite{KM76}.
To this end, in Section \ref{sec:detf}, we provide a detailed review of the category of lines with signs and the construction of the determinant line. 
This approach, however, may present additional challenges for those unfamiliar with the theory.

Alternatively, one could work within the category of ordinary lines.
In this case, both $ \sigma_{p}^{p-1} $ and $ s_{p} $ are only well-defined up to a factor $ \pm 1$.
The ambiguity disappears when we consider the partition function $ Z_{p} \left(s_{0 }, s^{L}_{D}, s^{E}_{D} \right) $.
However, motivated by recent developments in torsion theory (e.g., \cite{BismutShen24Arxiv}),  we believe that determining $ s_{p} $ itself is of intrinsic interest.
\end{re} 

\subsection{Proofs of (\ref{eq:zjj1}) and (\ref{eq:aaaz})}\label{isProofs}
The equations (\ref{eq:zjj1}) and (\ref{eq:aaaz}) are shown by the Quillen metric technique.

If $ \tau_{p} \left(\omega^{X}, h^{L}, h^{E} \right)$ is the holomorphic torsion of $ L^{p} \otimes E$ with respect to the metrics $ \left(\omega^{X}, h^{L}, h^{E} \right)$, then the Quillen metric $ \left\| \cdot \right\|^{{\rm Q}}_{\lambda_{p} \left(E\right)}  $ on $ \lambda_{p} \left(E\right)$ is defined by
\begin{align}\label{eq:tybw}
	\left\| \cdot \right\|^{{\rm Q}}_{\lambda_{p} \left(E\right)}= \exp\left(\tau_{p} \left(\omega^{X}, h^{L}, h^{E} \right)\right) \left| \cdot \right|^{}_{\lambda_{p} \left(E\right)}.
\end{align}
Thanks to \eqref{eq:tybw}, the proofs of  (\ref{eq:zjj1}) and (\ref{eq:aaaz}) can be reduced to show similar results for  $ 2\tau_{p} \left(\omega^{X}, h^{L}, h^{E} \right)$ and $ \log \left\| s_{p}  \right\|^{{\rm Q},2}_{\lambda_{p} \left(E\right)} $.

The full asymptotic expansion of $ 2\tau_{p} \left(\omega^{X}, h^{L}, h^{E} \right)$ was obtained by Finski \cite{Finski18}, following the earlier contribution of Bismut-Vasserot \cite{BismutVasserot89}.
In our case, the expansion is of type $ p^{i} \log p $ and $ p^{i} $ with $ i\in 1-\mathbf{N} $.
Under the prequantization condition, the term $ \left(\zeta'\left(-1\right)+\frac{\log\left(2\pi \right)}{12}+\frac{7}{24}\right){\rm rk}\left(E\right)\chi\left(X\right)+\frac{1}{2}\deg\left(E\right)$ in \eqref{eq:aaaz} is exactly from the constant term of this expansion. 

Thanks to our construction of $ s_{p} $, up to local normalisation, the evaluation of $ \log \left\| s_{p}  \right\|^{{\rm Q},2}_{\lambda_{p} \left(E\right)}$ 
is reduced to compute $ \log \left\| \sigma^{p-1}_{p}   \right\|^{{\rm Q},2}$, where $ \left\| \cdot \right\|^{{\rm Q},2} $ is defined by putting the both sides of \eqref{eq:3ko2} the Quillen metrics on $ \lambda_{p-1} \left(E\right), \lambda_{p} \left(E\right)$ and the obvious metric on $ \det \left(\left(L^{p} \otimes E\right)_{|D}\right) $.
The advantage of the Quillen metric over the $ L^{2} $-metric is that 
$ \log \left\| \sigma^{p-1}_{p}   \right\|^{{\rm Q},2}$ can be explicitly evaluated.

Assume now $ p= 1$.
If $ \left(\omega^{X}, h^{L}, h^{E} \right)$ satisfies Assumption A of Bismut \cite{B90immersion}, we can compute $ \log \left\| \sigma_{1}^{0}  \right\|^{{\rm Q},2}$ using Bismut-Lebeau's embedding formula \cite{BL91}.
For general metrics, we can reduce to the previous case by Bismut-Gillet-Soul\'e's anomaly formula \cite{BGS3}.
In this way, we establish an explicit formula for $ \log \left\| \sigma_{1}^{0}  \right\|^{{\rm Q},2}$  in Theorem \ref{thm:qbsp2}.

Ultimately, in Proposition \ref{prop:suof} and Theorem \ref{thm:ngzf},  we show that $  \log \left\| s_{p}   \right\|^{{\rm Q},2}_{\lambda_{p} \left(E\right)} $ is a polynomial in $ p$ of degree $ 2$.
Its constant term is given by $ \log \left\| s_{0 }   \right\|^{{\rm Q},2}_{\lambda_{0 } \left(E\right)} =  \log \left| s_{0 }   \right|^{2}_{\lambda_{0 } \left(E\right)} + 2\tau\left(\omega^{X}, h^{E} \right)$, which contributes to $ a_{0 } $.

In our approach, the geometrical Zabrodin-Wiegmann conjectured term arises naturally. 
However, from a purely $ L^{2} $-point of view, it seems to us that identifying such a term a priori would be highly challenging.

\subsection{Remarks on our assumption on the divisor $ D$ }\label{sReD}
In our paper, we restrict ourselves to the case where $ D$ is an effective reduced divisor.

If $ D$ is merely effective, it can be expressed as $ D= \sum_{i}^{} D_{i} $, a finite sum of effective reduced divisors. 
By iteratively applying the constructions and proofs described in Sections \ref{isConstr} and \ref{isProofs}, we obtain a modified version of Theorem \ref{thm:e2jm}. 
Replacing $ p$ by $ 2p$ in Theorem \ref{thm:e2jm} yields an example of such results.

If $ D$ is only assumed to be positive, then for sufficiently large $ p_{0 }\in \mathbf{N} $, $ p_{0 } D$ is equivalent to an effective divisor $ D^{\prime } $.
For each $0 \le r < p_{0 } $, applying the above results to $ D^{\prime } $ and $E^{\prime } =  L^{r} \otimes E$, we can obtain another version of Theorem \ref{thm:e2jm} for $ L^{p_{0 } k+r} \otimes E$ as $ k \to + \infty$.

\subsection{Comparison with related works}
One of the novelties of our work is the canonical construction, up to normalization data, of the element $ s_{p} \in \lambda_{p} \left(E\right)$ or equivalently the Slater determinant.
This construction is notably independent of any choice of metric data.

When $ g= 0 $ or $1$, in Sections \ref{sXg0} and \ref{sXg1}, we provide explicit expressions for our section $ s_{p} $, and show their compatibility with the constructions using bases of lowest Landau level eigenspaces, as given for instance in \cite[(2.20)]{Klevtsov16} for $ g= 0 $, and \cite[(14)]{AvronSeilerZograf95}, \cite[(2.34)]{Klevtsov16}, \cite[(39)]{Perice24} for $ g= 1$.
It is worth noting that our construction relies on another basis obtained via an inductive argument, which is different from the one used in the above references when $ g= 1$. 
This allows us to apply Bismut-Lebeau's embedding formula.
Our main result, Theorem \ref{thm:e2jm}, thus applies to these classical cases with general metrics, beyond the round metric when $ g= 0 $, or the flat metric when $ g= 1$.

In the absence of a canonical section $ s_{p}\in \lambda _{p} \left(E\right)$, in \cite{KlevtsovMaMarinescuWiegmann17}, the authors studied instead $\log \left(Z_{p } / Z_{p,0 }\right) $, where $ Z_{p,0 } $ is a reference partition function defined using fixed background metrics $ \left(\omega^{X}_{0 } , h^{L}_{0 } , h^{E}_{0}  \right)$.
When $ L$ is effective and reduced, our Theorem \ref{thm:e2jm} applies to $ \log \left(Z_{p} \left(s_{0 }, s_{D}^{L}, s_{D}^{E} \right)/Z_{p,0 } \left(s_{0 }, s_{D}^{L}, s_{D}^{E} \right)\right)$.
In particular, we recover \cite[Theorem 1]{KlevtsovMaMarinescuWiegmann17}.
For general positive line bundles $ L$, their results can also be deduced by incorporating the remarks in Section \ref{sReD}.

Our methods, however, do not extend to the case $ \beta \neq 1$, which remains the most challenging, as the interpretation of the partition function via the Quillen metric is no longer available. 
We leave this for future work.

Our results do not apply to the original Zabrodin-Wiegmann conjecture \cite{ZabrodinWiegmann06}, which concerns the non-compact Riemann surface. 
For further developments in that direction, we refer the reader to \cite{ByunKangSeo23,ByunSeoYang24} and \cite[Section 9.3]{Serfaty24}.

Finally, in the context of the free energy of the Coulomb gas, Bourgoin recently establishes a related result \cite{Bourgoin25}, analogous to ours but obtained through a different strategy.

\subsection{Organisation of the paper} 
This paper is organised as follows.
In Section \ref{S:Prel}, we introduce the necessary background on the determinant line of the cohomology and various metrics associated with it, including the Quillen metric and the $ L^{2} $-metric.

In Section \ref{Scans}, we construct the canonical element $ s_{p}\in  \lambda_{p} $.

In Section \ref{S:ZW}, we state our main result in a more general setting and provide its proof, leaving the evaluation of the Quillen norm of $ \sigma_{1}^{0 }$ to Section \ref{S:EQN}.

Finally, in Section \ref{S:EQN}, we evaluate the Quillen norm of $ \sigma_{1}^{0 }  $ and complete the proof of our main result.
 
In the whole paper, we use the superconnection formalism of Quillen \cite{Quillensuper}.
If $E = E^+ \oplus E^-$ is a $\mathbf{Z} _2$-graded vector space, the algebra $\End(E)$ is $\mathbf{Z} _2$-graded. If $\tau = \pm 1$ on $E^\pm$, if $a \in \End(E)$, the supertrace ${\rm Tr_{s} }[a]$ is defined by
\begin{align}
 {\rm Tr_ s}[a]={\rm Tr}[\tau a].
\end{align}
We adopt the convention that $\mathbf{N}=\{0,1, 2,\cdots\}$,  $ \mathbf{N} ^{*} =\{1, 2,\cdots\} $, $ \mathbf{R} _{+} = [0, +\infty)$, and $ \mathbf{R} _{+}^{*} = (0, +\infty)$.

\subsection{Acknowledgement}
The geometric Zabrodin-Wiegmann conjecture was introduced to us by Semyon Klevtsov during an ANR meeting in Orléans. We are deeply grateful to him for his enlightening discussions on quantum Hall effects and his valuable suggestions that have greatly improved this paper.
We are indebted to the referee for their insightful and constructive comments.
We acknowledge the partial financial support from NSFC under Grant No.~12371054.
S.S. is also partially supported by ANR Grant ANR-20-CE40-0017.

\section{Preliminary}\label{S:Prel}
The purpose of this section is to review some classical constructions related to the determinant line, the Ray-Singer holomorphic torsion, the Quillen metric, and the partition function associated to Slater determinants. 

This section is organised as follows.
In Section \ref{sec:detf}, we recall the determinant formalism developed by Grothendieck and Knudsen-Mumford \cite{KM76}.
 
In Section \ref{sec:hotq}, we review the classical theory of the Ray-Singer holomorphic torsion and Quillen metric associated to a Hermitian holomorphic vector bundle $ E$ on a Riemann surface $ X$.

Finally, in Section \ref{sec:pfunc}, we introduce the Slater determinant and the partition function for a family of smooth sections $ \left\{s_{i} \right\}_{1\le i \le k} $ of $ E$.

\subsection{Determinant functors}\label{sec:detf}
Let us review the determinant formalism following \cite[Chapter I]{KM76}.
Let $ \mathcal{L}_{{\rm is}} $ be the category of $ \mathbf{Z}_{2}$-graded complex lines and even isomorphisms.
Denote by $ \mathbf{C} $ the even trivial line in $ \mathcal{L}_{{\rm is}} $.

If $ L$ is an object in $ \mathcal{L}_{{\rm is}} $, denote by $ \epsilon \left(L\right)\in \mathbf{Z}_{2} $ its parity.
If $ L_{1}, L_{2} $ are two objects in $ \mathcal{L}_{{\rm is}} $, there is a well-defined $ \mathbf{Z}_{2} $-graded tensor product $ L_{1} \otimes L_{2} $  in $ \mathcal{L}_{{\rm is}} $.
The parity of $ L_{1} \otimes L_{2} $ is $ \epsilon\left(L_{1} \right) +\epsilon\left(L_{2}\right) $.
Moreover, $ L_{1} \otimes L_{2} $ is canonically identified with $ L_{2} \otimes {L_{1} } $ via 
\begin{align}\label{eq:w2ny}
 \ell_{1} \otimes \ell_{2} \in L_{1} \otimes L_{2} \to  \left(-1\right)^{\epsilon\left(L_{1}\right) \epsilon\left(L_{2}\right) } \ell_{2} \otimes \ell_{1} \in L_{2} \otimes L_{1}. 
\end{align}

Let $ L^{-1} $ be the right inverse of $ L$.
There is a canonical isomorphism
\begin{align}\label{eq:qxd2}
 L\otimes L^{-1} \simeq \mathbf{C}. 
\end{align}
If $ s\in L$, then $ \left(s\right)_{{\rm r}} ^{-1}\in L^{-1} $ is the unique right inverse such that 
\begin{align}\label{eq:glaf}
	s \left(s\right)_{{\rm r}} ^{-1}= 1.
\end{align}
Put $ \left(s\right)_{{\rm l}} ^{-1}= \left(-1\right)^{\epsilon \left(L\right)} \left(s\right)_{{\rm r}} ^{-1}. $
By \eqref{eq:w2ny} and \eqref{eq:glaf},
\begin{align}\label{eq:trqp}
 \left(s\right)_{{\rm l}}^{-1} s = 1.
\end{align}

Let $ \mathcal{V}_{\rm is}  $ be the category of finite dimensional complex vector spaces and isomorphisms. 
If $ n\in \mathbf{Z} $, denote $ [n]_{2} $ the corresponding class in  $ \mathbf{Z}_{2} $.
If $ V$ is an object in $ \mathcal{V}_{{\rm is}} $, put 
\begin{align}\label{eq:omds}
	\det\left(V\right)= \left( \Lambda^{\dim V} \left(V\right),[\dim V]_{2} \right).
\end{align}
Then, $ \det$ induces a functor from $ \mathcal{V}_{{\rm is}} $ to $ \mathcal{L}_{{\rm is}} $.

If $ L$ is an ordinary line, denote by $ L_{{\rm ev}} $ and $ L_{{\rm odd}} $ the corresponding even and odd lines in the category $ \mathcal{L}_{{\rm is}} $.
By \eqref{eq:omds}, we have 
\begin{align}\label{eq:sn3z}
 \det \left(L \otimes V\right)=\left(L^{\dim V} \otimes \Lambda^{\dim V} \left(V\right),\left[\dim V\right]_{2} \right)= \left(L_{{\rm ev}}\right)^{\dim V}   \otimes \det\left(V\right). 
\end{align}
When $ V= \mathbf{C} $, the above formula reduces to 
\begin{align}\label{eq:oywi}
 \det \left(L\right)= L_{{\rm ev}} \otimes \mathbf{C}_{{\rm odd}}  = L_{{\rm odd}}.
\end{align}


Let $ \mathcal{C}_{{\rm qi}} $ be the category of bounded complexes of finite dimensional vector spaces and quasi-isomorphisms.
By \cite[Theorem 1]{KM76}, the determinant functor extends to $ \mathcal{C}_{{\rm qi}} $, 
\begin{align}\label{eq:oe3g}
 \det : \mathcal{C}_{{\rm qi}} \to \mathcal{L}_{{\rm is}}. 
\end{align}
If $ E$ is an object in $ \mathcal{C}_{{\rm qi}} $, if $ p,q \in \mathbf{Z} $ with $ p\le q$ so that $ E= \oplus_{i= p}^{q} E^{i} $, then 
\begin{align}\label{eq:vgcd}
 \det\left(E\right)= \bigotimes_{i= p }^{q}  \det\left(E^{i} \right)^{\left(-1\right)^{i}},
\end{align}
where the tensor product is taken in order from $ p$ to $ q$.

If $E$ is acyclic, the fundamental properties of determinant functors \cite[p.~25]{KM76} imply $ \det E \simeq \mathbf{C} $.
Equivalently, $ \det E$ is equipped with a canonical non-zero element $ \tau$.
If $ \delta $ is the differential on $ E$, if $ i\in \mathbf{Z} $, and if $ e^{i} $ is a non-zero element in $ \det \left(E^{i} /\delta E^{i-1}\right) $, then 
\begin{align}\label{eq:lwdj}
 \tau= \bigotimes  _{i= p}^{q}  \left(\delta e^{i-1} \wedge e^{i}\right)^{\left(-1\right)^{i} }.
\end{align}
In the above formula, if $ i$ is odd, we use the convention
\begin{align}\label{eq:pgqy}
	\left(\delta e^{i-1} \wedge e^{i}\right)^{-1 }= \left(\delta e^{i-1}\right)^{-1}_{{\rm r}} \wedge \left(e^{i} \right)^{-1}_{{\rm l}}. 
\end{align}
Clearly, the right hand side of \eqref{eq:lwdj} is independent of the choice of $ e^{i}$.

\begin{re}\label{re:gkjm}
 The category $ \mathcal{L}_{{\rm is}} $ is introduced by Knudsen-Mumford \cite{KM76} to resolve sign inconsistencies when studying the compatibility of $ {\rm det}$ with direct sum of vector spaces.
 The precise determinant functor is a pair $ \left({\rm det},{\rm i}\right)$, where $ {\rm i}$ gives the desired compatibility. 
 For further details, we refer readers to \cite{KM76}.

 Alternatively, one can work in the category of ordinary lines and isomorphisms.
 In this case, the section $ \tau$ in \eqref{eq:lwdj} is only well-defined up to multiply by $ \pm1$ (see \cite[Remark 1.2]{BGS1}).
\end{re}

\subsection{Holomorphic torsion and Quillen metric}\label{sec:hotq}
Let $ X$ be a closed connected Riemann surface of genus $ g$.
Denote by $ TX $ and $ \overline{TX} $ the holomorphic and anti-holomorphic tangent bundles of $ X$.

Let $ H_{{\rm dR}}  \left(X,\mathbf{C} \right)$ be the de Rham cohomology of $ X$.
Since $ X$ is K\"ahler, $ H_{{\rm dR}}  \left(X,\mathbf{C} \right)$ is bigraded so that     
\begin{align}\label{eq:qawr}
&\dim H^{0, 0 }_{{\rm dR}}  \left(X,\mathbf{C} \right)= 1,&\dim H^{0, 1}_{{\rm dR}}  \left(X,\mathbf{C} \right)= g,\\
&\dim H^{1, 0 }_{{\rm dR}}  \left(X,\mathbf{C} \right)= g,&\dim H^{1,1}_{{\rm dR}}  \left(X,\mathbf{C} \right)= 1.\notag
\end{align}

Let $ E$ be a holomorphic vector bundle on $ X$.
Let $ {\rm rk}\left(E\right) \in \mathbf{N} $ be the rank of $ E$.
Let $ {\rm deg}\left(E\right)\in \mathbf{Z}$ be the degree of $ E$.
If  $ c_{1} \left(E\right)\in H^{1,1}_{{\rm dR}}  \left(X,\mathbf{C} \right)$ is the first Chern class of $ E$, then
\begin{align}\label{eq:ipog}
	{\rm deg}\left({E}\right) = \int_{X} c_{1} \left(E\right).
\end{align}
 

Let $ \mathcal{O}_{X} $ be the sheave of holomorphic functions on $ X$.
Let $ \mathcal{O}_{X} \left(E\right)$ be the sheave of holomorphic sections of $ E$.
Let $ H \left(X,\mathcal{O}_{X}  \left(E\right)\right)= H^{0 }  \left(X,\mathcal{O}_{X}  \left(E\right)\right)\oplus H^{1} \left(X,\mathcal{O}_{X}  \left(E\right)\right)$ be the cohomology of $ \mathcal{O}_{X} \left(E\right)$.
When there is no risk of confusion, we will simply write $ H \left(X,E\right)$ for $ H \left(X,\mathcal{O}_{X} \left(E\right)\right)$.


Let $ \chi\left(E\right) \in \mathbf{Z}  $ be the holomorphic Euler characteristic of $ E$, i.e., 
\begin{align}\label{eq:q3ka}
	\chi\left(E\right)= \dim H^{0  } \left(X,E \right)-\dim H^{ 1 } \left(X,E \right).
\end{align}
By Riemann-Roch theorem, we have
\begin{align}\label{eq:gbhd1}
	\chi\left(E\right) =  \left(1-g\right){\rm rk}\left({E}\right)  +{\rm deg}\left({E}\right).
\end{align}

Set  
\begin{align}\label{eq:evqq}
 \lambda \left(E\right)= \det H^{0 } \left(X,E\right) \otimes \left(\det H^{1} \left(X,E\right)\right)^{-1}.  
\end{align}
Then, $ \lambda \left(E\right)$ is an object in $ \mathcal{L}_{{\rm is}} $.
The parity of $ \lambda \left(E\right)$ is given by $ [\chi\left(E\right) ]_{2} \in \mathbf{Z}_{2} .$ 

Let $ \left(\Omega^{0, \bullet} \left(X,E\right),\overline{\partial}^{E} \right)$ be the Dolbeault complex with values in $ E$.
Classically,  the cohomology of $ \left(\Omega^{0, \bullet} \left(X, E\right),\overline{\partial}^{E} \right)$ coincides with $ H^{\bullet} \left(X,E\right)$.

Let $ \omega^{X} $ be a K\"ahler form on $ X$. 
It induces a volume form on $ X$, 
\begin{align}\label{eq:epzu}
 dv_{X} = \frac{\omega^{X}}{2\pi}. 
\end{align}
Let $ h^{E} $ be a Hermitian metric on $ E$.
Then, $ \left(\omega^{X}, h^{E}\right) $ induces a Hermitian metric $ \left\langle,\right\rangle _{\Lambda^{\bullet} \left(\overline{T^{*} X} \right)\otimes {E}}$ on $ \Lambda^{\bullet} \left(\overline{T^{*} X} \right)\otimes {E}$. 
Ultimately, we get an $ L^{2} $ metric on $ \Omega^{0, \bullet} \left(X,E\right)$ defined by 
\begin{align}\label{eq:si2u}
\left\langle s_{1} , s_{2} \right\rangle _{L^{2} } = \int_{z\in X} \left\langle s_{1}\left(z\right) , s_{2}\left(z\right) \right\rangle _{\Lambda^{\bullet} \left(\overline{T^{*} X} \right)\otimes E} dv_{X}.
\end{align}
The corresponding norm will be denoted by $ |\cdot |_{L^{2} }$.

Let $ \overline{\partial}^{E*} $ be the adjoint of $ \overline{\partial}^{E} $ with respect to the above $ L^{2} $ metric. 
Put
\begin{align}\label{eq:ztwh}
	\Box^{E} = \overline{\partial}^{E*}\overline{\partial}^{E}+\overline{\partial}^{E}\overline{\partial}^{E*}.
 \end{align} 
Then, $ \Box^{E} $ is the Kodaira Laplacian.
By Hodge theory, we have
\begin{align}\label{eq:nigg}
 \ker \Box^{E} \simeq H \left(X,E\right).
\end{align}

For $ s\in \mathbf{C} $ and $ {\rm Re}\,s> 1$, set
\begin{align}\label{eq:h3zm}
 \zeta^{E}  \left(s\right)= \sum_{\lambda \in \mathbf{R}_{+}^{*} \cap \Sp \Box^{E} _{|\Omega^{0,1 } \left(X,E\right)}  } \frac{1}{\lambda^{s} },
\end{align}
where the sum is taken over the positive spectrum of $ \Box^{E} _{|\Omega^{0,1 } \left(X,E\right)}$, counted with multiplicity.
By a result of \cite{Seeley66}, $ \zeta^{E}  \left(s\right)$ has a meromorphic extension to $ s\in \mathbf{C} $ and is holomorphic at $ s= 0 $.

\begin{defin}\label{def:vsaz}
 The Ray-Singer holomorphic torsion \cite{RScomplextor} of $ E$ with respect to $ \left(\omega^{X}, h^{E} \right)$  is defined by 
 \begin{align}\label{eq:2vdr}
	\tau\left(\omega^{X},h^{E}\right)= \frac{1}{2}\frac{\partial \zeta^{E} }{\partial s}\left(0 \right)  \in \mathbf{R} .
 \end{align} 
\end{defin}

By \eqref{eq:evqq} and \eqref{eq:nigg}, $ \lambda \left(E\right)$ equips an $ L^{2} $-metric $ |\cdot|_{\lambda \left(E\right)} $.

\begin{defin}\label{def:yogw}
 The Quillen metric \cite{Quillendeter,BGS3} of $ \lambda \left(E\right)$ is defined by
 \begin{align}\label{eq:iybw}
	\left\| \cdot \right\|^{{\rm Q}}_{ \lambda \left(E\right)}  =  |\cdot|_{\lambda \left(E\right)}  \exp\left({\tau\left(\omega^{X},h^{E}\right)}\right). 
 \end{align}
\end{defin}

\subsection{Slater determinant and partition function}\label{sec:pfunc}
If $ k\in \mathbf{N}^{*} $, let $ \mathfrak{S}_{k}  $ be the permutation group of order $ k$. If $ \sigma\in \mathfrak{S}_{k} $, denote by $ \epsilon \left(\sigma \right) \in \{\pm1\}$ the sign of $ \sigma $.

Let $ s= \left(s_{i}\right)_{1\le i\le k} $ be a family of smooth sections of $ E$.
If $ z=  \left(z_{i} \right)_{1\le i\le k} \in  X^{k}$, denote 
\begin{align}\label{eq:wmbo}
 \det\left(s \left(z \right)\right)= \sum_{\sigma \in \mathfrak{S} _{k} } \epsilon \left(\sigma \right)\prod_{i= 1}^{k} s_{\sigma{i} } \left(z_{i} \right) \in \bigotimes_{i= 1}^{k} E_{z_{i} } . 
\end{align}
Then, $ \det\left(s\right)$ is a section of $ E^{\boxtimes k} $ over $ X^{k} $.

\begin{defin}\label{def:akdc}
 The section $ \frac{1}{\sqrt{k!}} \det\left(s\right)\in C^{\infty} \left(X^{k},  E^{\boxtimes k}\right)$ is called the normalised Slater  determinant of $ s= \left(s_{i} \right)_{1\le i \le k} $.
\end{defin}

\begin{exa}\label{exa:akdc}
 Let $ u$ be a smooth section of $ E$.
 If $ \left(f_{i} \right)_{1\le i \le k} $ is a family of smooth functions on $ X$, then the Slater determinant of $ \left( f_{i} u\right)_{1\le i \le k}$ is given by 
 \begin{align}\label{eq:lzei}
	\det \left[\left(f_{j} \left(z_{i} \right)\right)_{1\le i,j\le k} \right] \cdot \prod_{i= 1}^{k}  u\left(z_{i} \right).
 \end{align}
\end{exa}

Clearly, $ h^{E} $ induces a Hermitian metric $ h^{E^{\boxtimes k} } $ on $ E^{\boxtimes k} $.
The corresponding norm will be denoted by $ |\cdot |_{h^{E^{\boxtimes k} } } $.
Let $ dv_{X^{k} } $ be the volume form on $ X^{k} $ induced by $ dv_{X} $.

\begin{defin}\label{def:ye12}
 Set
 \begin{align}\label{eq:owsl}
	Z\left(s,\omega^{X}, h^{E} \right)= \frac{1}{k!}\int_{z \in X^{k} } |\det\left(s\left(z\right)\right)|^{2}_{h^{E^{\boxtimes k} } } dv_{X^{k} }. 
 \end{align}
\end{defin}

Let $ S_{k} $ be a symmetric $ k \times k$ matrix, where $ \left(i,j\right)$-entry is given by  $ \left\langle s_{i} , s_{j} \right\rangle _{L^{2} } $.
For sake of completeness, we give a proof of the following well-known result \cite[(2.6)]{Klevtsov14}.

\begin{prop}\label{prop:apm1}
 The following identity holds, 
 \begin{align}\label{eq:xm1t}
	Z\left(s,\omega^{X}, h^{E} \right)= \det\left(S_{k}\right)= \left| s_{1} \wedge\ldots \wedge s_{k}  \right|^{2}_{L^{2}}.
 \end{align}
\end{prop}
\begin{proof}
 By \eqref{eq:wmbo}, for $ z= \left(z_{i} \right)_{1\le i\le k}\in X^{k} $, we have
 \begin{align}\label{eq:dnn1}
	\left|\det\left(s \left(z \right)\right)\right|^{2}_{h^{E^{\boxtimes k}} } = \sum_{\sigma,\sigma' \in \mathfrak{S} _{k} } \epsilon\left(\sigma \right)\epsilon\left(\sigma' \right)\prod_{i= 1}^{k} \left\langle s_{\sigma{i} } \left(z_{i} \right), s_{\sigma'{i} } \left(z_{i} \right)\right\rangle . 
 \end{align}
 By \eqref{eq:dnn1} and by Fubini's theorem, we get 
 \begin{align}\label{eq:gaoc}
\int_{X^{k} }^{}	\left|\det\left(s \left(z \right)\right)\right|^{2}_{h^{E^{\boxtimes k} } } dv_{X^{k} } = \sum_{\sigma,\sigma' \in \mathfrak{S} _{k} } \epsilon\left(\sigma \right)\epsilon\left(\sigma' \right)\prod_{i= 1}^{k}  \left\langle s_{\sigma{i} } , s_{\sigma'{i} } \right\rangle _{L^{2} }. 
 \end{align}
 From \eqref{eq:owsl} and \eqref{eq:gaoc}, we get the first identity in \eqref{eq:xm1t}.

 The second identity in \eqref{eq:xm1t} is a trivial identity in linear algebra.
 The proof of our proposition is complete.
\end{proof}

\section{The canonical section in the determinant line}\label{Scans}
Given an effective reduced divisor $ D$ on $ X$, we can associate a positive holomorphic line bundle $ L$ and a canonical holomorphic section $ s_{D} $ of $ L$.
The purpose of this section is to construct, for each $ p \ge 1$, a non-zero element $ s_{p} $ in $ \det H(X,L^{p} \otimes E)$, using $ s_{D} $, an initial data $ s_{0 } \in \det H(X,E)$, and some local normalisation data defined over $ D$.
We emphasise that the construction does not rely on any metric data.

This section is organised as follows.
In Section \ref{sPLB}, we recall some well-known results on positive line bundles.

In Section \ref{sDiv}, we review the basic properties of effective reduced divisors, the associated positive line bundle, and the canonical holomorphic section.

In Section \ref{sec:1c1y}, we construct our canonical non-zero section $ s_{p} \in \det H\left(X,L^{p} \otimes E\right)$.

In Section \ref{sEtri}, we specialize this construction to the case where $ E= \mathbf{C} $ is trivial.
In this case, up to multiply by a complex number of module $ 1$, the initial data $ s_{0 }\in \det H\left(X,\mathbf{C} \right)$ is canonically determined by some $ \mathbf{Z} $-structure from the Hodge theory.

Finally, in Sections \ref{sXg0} and \ref{sXg1}, we explicitly evaluate $ s_ p$ for the cases $ g= 0$ and $ g= 1$, respectively.

\subsection{A positive line bundle}\label{sPLB}
Let $ L$ be a positive line bundle on $ X$.
By \cite[p.~214]{GH94}, this is equivalent\footnote{Since $ X$ is K\"ahler, the $ \partial \overline{\partial} $-lemma holds, so that $ L$  is positive if and only if $ c_{1} \left(L\right)\in H^{1,1}_{{\rm dR}}  \left(X,\mathbf{R} \right)$  is positive. 
Since $ X$ is connected and has complex dimension $ 1$, the integration over $ X$  identifies the lines $ H^{1,1}_{{\rm dR}}  \left(X,\mathbf{R} \right) \simeq \mathbf{R} $ together with their positivity. Therefore, the positivity of $ c_{1} \left(L\right)$  coincides with the positivity of $ {\rm deg}\left(L\right)$.
} to  
\begin{align}\label{eq:2yxi}
 {\rm deg}\left(L\right)> 0. 
\end{align}

Let $ p\in \mathbf{Z} $.
We consider the constructions from Section \ref{S:Prel}, replacing $E$ with $L^p \otimes E$ and using corresponding notation.
In particular, write
\begin{align}\label{eq:xfpn}
\lambda_{p}\left(E\right)= \lambda\left(L^{p} \otimes E\right).
\end{align}

Since $ L$ is positive, by the vanishing theorem \cite[Theorem B, p.~159]{GH94},  there exists $ p_{0 }\in \mathbf{Z} $ such that if $ p\ge p_{0 } $,  
\begin{align}\label{eq:rpfs}
 H^{1}\left(X,L^{p} \otimes E\right)= 0.
\end{align}
 
By \eqref{eq:gbhd1} and \eqref{eq:rpfs}, for $ p\ge p_{0 } $,
\begin{align}\label{eq:gbhd2}
	\dim H^{0}\left(X,L^{p} \otimes E\right) = \left(1-g +p \deg\left(L\right) \right){\rm rk}\left(E\right) + {\rm deg}\left(E\right). 
\end{align}
Moreover, 
\begin{align}\label{eq:jcta}
 \lambda_{p} \left(E\right)= \det H^{0}\left(X,L^{p} \otimes E\right).
\end{align}

\begin{re}\label{re:bxyk}
 If $ a\in \mathbf{R} $, denote $ [a]$ the integer part of $ a$.
 If $ E= \mathbf{C} $ is the trivial vector bundle, by \cite[p.~214]{GH94}, one can take $ p_{0 } = \left[\left(2g-2\right)/|D|\right]+1$.
\end{re}

\subsection{An effective reduced divisor}\label{sDiv}
Let $ D\subset X$ be a non empty finite subset of $ X$.
We identify $ D$ with an effective reduced divisor $ \sum_{z \in D} z$ on $ X$.

Let $L$ be the holomorphic line bundle on $ X$ associated to $ D$ \cite[p.~134]{GH94}. 
By \cite[p.~214]{GH94}, 
\begin{align}\label{eq:wwhc}
 \deg L = |D|.
\end{align}
By \eqref{eq:2yxi} and \eqref{eq:wwhc}, $ L$ is positive.

If $ s$ is a non-zero meromorphic section of $ L$, denote by $ {\rm div}\left(s\right)$ the associated divisor.
By \cite[p.~136]{GH94}, up to multiplication by a non-zero complex number, there is a unique  holomorphic section $ s_{D} $ of $ L$ such that 
\begin{align}\label{eq:mezs}
 {\rm div}\left(s_{D}\right)= D.
\end{align}

Let $ \mathcal{M}_{pD}\left(E\right) $ be the space of meromorphic sections of $ E$ with vanishing order $ \ge -p$ at each $ z\in D$.
The same reference above gives an isomorphism of vector spaces, 
\begin{align}\label{eq:xg2h}
 f\in \mathcal{M}_{pD} \left(E\right)\to s_{D}^{p} f\in H^{0 } \left(X,L^{p} \otimes E\right).
\end{align}

\begin{re}\label{re:tvgq}
 Most of the above results hold for general divisors.
 The effectiveness of $ D$ induces a filtration on the space of meromorphic sections of $ E$ with possible poles along $ D$, 
 \begin{align}\label{eq:alri}
	 \cdots \subset \mathcal{M} _{\left(p-1\right)D} \left(E\right) \subset  \mathcal{M}_{pD} \left(E\right) \subset \cdots 
 \end{align}
 This filtration is used essentially but implicitly in our construction of the canonical section in $ \lambda_{p} \left(E\right)$.
\end{re}

In the sequel, we will fix such an effective reduced divisor $ D$, together with the associated line bundle $ L$ and the canonical holomorphic section $ s_{D} $.
We consider only the non negative tensor power of $ L$, i.e., $ p \in \mathbf{N} $.

\subsection{The canonical section in $ \lambda_{p} \left(E\right)$ }\label{sec:1c1y}
The holomorphic vector bundle $ L^{p} \otimes E$ restricts to $ D$.
We identity this restriction $ \left(L^{p} \otimes E\right)_{|D} $ with the obvious sheaf on $ D$.

Let $ \iota : D\to X $ be the natural embedding.
The direct image $ \iota_{*}\left(\left(L^{p} \otimes E\right)_{|D}\right) $ is a sheaf on $ X$ supported on $ D$, whose fibre over $z\in D$ is given by $ L^{p}_{z}  \otimes E_{z} $.

We have an exact sequence of $ \mathcal{O}_{X} $-modules,  
\begin{equation}\label{dia:pgtl}
	\begin{tikzcd}
		0 & \mathcal{O}_{X} \left(L^{p-1} \otimes E\right) &	\mathcal{O}_{X} \left(L^{p} \otimes E\right) & \iota_{*}\left(\left(L^{p} \otimes E\right)_{|D}\right) & 0,
		 \arrow["", from=1-1, to=1-2]
		 \arrow["s_{D} ", from=1-2, to=1-3]
		 \arrow["r_{D} ", from=1-3, to=1-4]
		 \arrow["", from=1-4, to=1-5]
	\end{tikzcd} 
\end{equation}
where the first morphism is given by the multiplication by $ s_{D} $ and the second one is the restriction to $ D$.

By \eqref{dia:pgtl}, we get the associate long exact sequence of vector spaces, 
\begin{equation}\label{dia:vooe}
	\begin{tikzcd}
		 \cdots & H^{\bullet} \left(X,L^{p-1} \otimes E\right) &	H^{\bullet} \left(X,L^{p} \otimes E\right) & H^{\bullet} \left(X,\iota_{*}\left(\left(L^{p} \otimes E\right)_{|D}\right) \right) & \cdots
		 \arrow["", from=1-1, to=1-2]
		 \arrow["", from=1-2, to=1-3]
		 \arrow["", from=1-3, to=1-4]
		 \arrow["", from=1-4, to=1-5]
	\end{tikzcd} 
\end{equation}
By \eqref{eq:lwdj}, there is a canonical non-zero section in the determinant line of the acyclic complex \eqref{dia:vooe}.
It induces therefore a canonical isomorphism $ \sigma^{p-1}_{p} $ in $ \mathcal{L}_{{\rm is}} $,
\begin{align}\label{eq:x1mn}
	\lambda_{p-1} \left(E\right)\otimes \det\left(H \left(X,\iota_{*}\left(\left(L^{p} \otimes E\right)_{|D}\right) \right)\right) \simeq \lambda_{p} \left(E\right).
\end{align}

Put 
\begin{align}\label{eq:iga2}
 \lambda_{{\rm ev}} \left( L_{|D}\right) = \bigotimes _{z \in D} L_{z,{\rm ev}}.
\end{align}
Then, $ \lambda_{{\rm ev}} \left( L_{|D}\right)$ is an even line.
 
Let $ \lambda \left(E_{|D}\right)$ be the determinant line of the vector bundles $ E_{|D} $ on $ D$.
When an order is fixed on the set $ D$, we have canonical isomorphism in $ \mathcal{L}_{{\rm is}} $,  
\begin{align}\label{eq:czyn}
	\lambda \left( E_{|D}\right) \simeq \bigotimes _{z \in D} \det\left(E_{z}\right). 
\end{align}

Observes that 
\begin{align}\label{eq:m1ux}
&	H^{0 } \left(X,\iota_{*}\left(\left(L^{p} \otimes E\right)_{|D}\right) \right)= \bigoplus_{z \in  D} L^{p}_{z}  \otimes E_{z}, &H^{1 } \left(X,\iota_{*}\left(\left(L^{p} \otimes E\right)_{|D}\right)\right)= 0. 
\end{align}
By \eqref{eq:sn3z}, \eqref{eq:iga2}-\eqref{eq:m1ux}, we get 
\begin{align}\label{eq:byiz}
 \det H\left(X,\iota_{*}\left(\left(L^{p} \otimes E\right)_{|D}\right) \right) \simeq \left(\lambda_{{\rm ev}}   \left( L_{|D}\right) \right)^{p{\rm rk}E} \otimes \lambda \left(E_{|D} \right).
\end{align}

\begin{re}\label{re:dsep}
Unlike \eqref{eq:czyn}, the isomorphism \eqref{eq:byiz} does not depend on the order of the points in $ D$. 
\end{re}

Using \eqref{eq:byiz}, we can rewrite the isomorphism $ \sigma_{p}^{p-1}$ in \eqref{eq:x1mn} as 
\begin{align}\label{eq:sqh3}
	 \lambda_{p-1} \left(E\right)\otimes \left(\lambda_{{\rm ev}} \left( L_{|D}\right) \right)^{p{\rm rk}E} \otimes \lambda \left(E_{|D} \right) \simeq \lambda_{p} \left(E\right).
\end{align}

Set
\begin{align}\label{eq:gunc}
 \sigma^{0}_{p} =  \prod_{i= p}^{1}   \sigma_{i}^{i-1}, 
\end{align}
where the product is taken in the inverse order from $ i= p$ to $ 1$.
By \eqref{eq:sqh3} and \eqref{eq:gunc}, $ \sigma^{0}_{p} $ gives a canonical isomorphism in $ \mathcal{L}_{{\rm is}} $, 
\begin{align}\label{eq:av1q}
 \lambda \left(E\right)\otimes \left(\lambda_{{\rm ev}} \left( L_{|D}\right) \right)^{p\left(p+1\right) {\rm rk}E/2} \otimes\left( \lambda \left(E_{|D} \right)\right)^{p}  \simeq \lambda_{p} \left(E\right).
\end{align}

\begin{defin}\label{def:zwfd}
	If $ s_{0 } \in \lambda \left(E\right)$, $ s_{D}^{L}    \in \lambda_{{\rm ev}} \left( L_{|D}\right)$, and $ s_{D}^{E}  \in \lambda \left(E_{|D} \right) $, set 
	\begin{align}\label{eq:bihv}
	 s_{p} =  \sigma_{p}^{0} \left(s_{0 } \otimes \left(s_{D}^{L}  \right) ^{p\left(p+1\right) {\rm rk}E/2}\otimes \left(s_{D}^{E} \right)^{p}   \right)\in  \lambda_{p} \left(E\right).
	\end{align}	 
\end{defin}

If $ p\ge p_{0 } $, by \eqref{eq:jcta}, $ s_{p} $ is an element in $ \det H^{0 } \left(X,L^{p} \otimes E\right)$.

In the sequel, we will always take non-zero data $ s_{0 } $, $ s_{D}^{L} $, and $ s_{D}^{E}$, so that $ s_{p} $ is also non-zero.

\begin{re}\label{re:bnhq}
Note that our construction of $ s_{p} $ does not involve any metric data.
\end{re}

\begin{re}\label{re:awko}
	In the context of real manifolds, a similar construction is used in \cite[Section 2.4]{ShenYuMorseSmale}, where the Smale filtration serves as an analogue to the filtration \eqref{eq:alri} here.
\end{re}

\subsection{The case where $ E$ is trivial}\label{sEtri}
In this section, we assume $ E$ is trivial.
Write  
\begin{align}\label{eq:hwjc}
 \lambda_{p} = \lambda \left(L^{p} \right).
\end{align}

Let $ H_{{\rm dR}}  \left(X,\mathbf{Z} \right)$ be the cohomology of the locally constant sheave $ \mathbf{Z} $.
Then, 
\begin{align}\label{eq:bxd2}
&	H^{0 }_{{\rm dR}}  \left(X,\mathbf{Z} \right)= \mathbf{Z}, &H^{1 }_{{\rm dR}}  \left(X,\mathbf{Z} \right) \simeq \mathbf{Z}^{2g}, && H^{2}_{{\rm dR}}  \left(X,\mathbf{Z} \right)= \mathbf{Z}.
\end{align}
Since $ X$ is connected and oriented, the first and last identities in \eqref{eq:bxd2} are canonical, while the second one is not, as notation indicates.

Similarly identities hold when we replace $ \mathbf{Z} $ by $ \mathbf{R} $.
We have the  canonical isomorphism of real vector spaces, 
\begin{align}\label{eq:qdaj}
	H^{\bullet}_{{\rm dR}}  \left(X,\mathbf{Z} \right) \otimes_{\mathbf{Z} }\mathbf{R} =   H^{\bullet}_{{\rm dR}}  \left(X,\mathbf{R} \right).
\end{align}
In particular, $ H^{\bullet }_{{\rm dR}}  \left(X,\mathbf{Z} \right)$ is a lattice in $ H^{\bullet }_{{\rm dR}}  \left(X,\mathbf{R} \right)$.

By Hodge theory \cite[VI.(11.3)]{DemaillyBook}, we have the isomorphism of complex vector spaces,\footnote{Recall that $ H^{1} \left(X,\mathbf{C} \right)$ is the cohomology of $ \mathcal{O}_{X} $.}
\begin{align}\label{eq:nztx}
 H^{1}_{{\rm dR}}  \left(X,\mathbf{R} \right)\otimes_{\mathbf{R} } \bC = H^{1} \left(X,\mathbf{C} \right)\oplus \overline{H^{1} \left(X,\mathbf{C} \right)}.
\end{align}
Then,  $ H^{1}_{{\rm dR}}  \left(X,\mathbf{R} \right)$ is equipped with a complex structure and hence is  canonically oriented.
 
An oriented generator of $H^{1}_{{\rm dR}}  \left(X,\mathbf{Z} \right)$ is  an oriented basis of $ H^{1}_{{\rm dR}}  \left(X,\mathbf{R} \right)$.
It defines a canonical element 
\begin{align}\label{eq:hqq3}
 \tau_{\mathbf{Z} }\in  \det H^{1}_{{\rm dR}}  \left(X,\mathbf{R} \right).
\end{align}

By \eqref{eq:nztx}, we have the canonical isomorphism in $ \mathcal{L}_{{\rm is}} $, 
\begin{align}\label{eq:gpty}
 \left(\det H^{1}_{{\rm dR}}  \left(X,\mathbf{R} \right)\right) \otimes_\mathbf{R} \mathbf{C} = \det H^{1} \left(X,\mathbf{C} \right) \otimes \overline {\det H^{1} \left(X,\mathbf{C} \right)}.
\end{align}
If $ s  \in \det H^{1} \left(X,\mathbf{C} \right)$, then  $ \sqrt{-1}  s  \otimes  \overline{s}  $ is a positive section in $ \det H^{1}_{{\rm dR}}  \left(X,\mathbf{R} \right)$.
Up to multiple by a complex number of module $ 1$, there is a unique $ s_{0 }^{1}  \in \det H^{1} \left(X,\mathbf{C} \right)$ such that 
\begin{align}\label{eq:dlig}
\sqrt{-1} 	s_{0 }^{1}  \otimes \overline{s}_{0 }^{1}  = \tau_{\mathbf{Z} }. 
\end{align}

Let $ s_{0 }^{0} $ be the positive generator of $ H^{0 }_{{\rm dR}} \left(X,\mathbf{Z} \right)$. 
We identify $ s_{0 }^{0} $ with the constant function $ 1$  on $ X$, which is an element in $ H^{0 } \left(X,\mathbf{C} \right)$.
Set
\begin{align}\label{eq:agmh}
 s_{0 } = s_{0 }^{0} \otimes \left(s_{0 }^{1} \right)^{-1}  \in \det H \left(X,\mathbf{C} \right).
\end{align}
With this choice of $ s_{0 }$, up to multiple by a complex number of module $ 1$, 
\begin{align}\label{eq:ehvt}
 s_{p} = \sigma^{0}_{p} \left(s_{0 } \otimes \left(s^{L}_{D}\right)^{p\left(p+1\right)/{2} } \right) \in  \lambda_{p}
\end{align}
depends only on the normalisation data $ s^{L}_{D} \in \lambda \left( L_{|D} \right)$.

\subsection{The case where $ X$ is $ \mathbf{CP}^{1}  $}\label{sXg0}
In this section, we assume that $ X= \mathbf{CP}^{1}$.
Then, $ g= 0 $.
If $ \left(x,y\right)\in \mathbf{C}^{2} \backslash \left\{\left(0, 0 \right)\right\}$, denote $ z = [x:y]\in X$ the corresponding point in homogeneous coordinates.

Let $ D$ be the reduced divisor defined by the point $ [0 : 1 ] \in  X$.
Then, $ L= \mathcal{O} \left(1\right)$. 
If $ p \ge 0 $, 
	\begin{align}\label{eq:tqfx}
&\dim H^{0 } \left(X,L^{p} \right)= p+1,&	 \dim H^{1} \left(X,L^{p} \right)= 0. 
	\end{align}
We can identify $ H^{0 } \left(X,L^{p} \right)$ with the space of homogenous polynomials on $ \left(x, y \right)\in \mathbf{C}^{2} $ of degree $ p$.
The canonical section $ s_{D} $ is defined by the polynomial $x$.

Moreover, $ \mathcal{M}_{pD} $ is the space of polynomial on $ y/x $ of degree $ \le p$.
The identification \eqref{eq:xg2h} is given by 
	\begin{align}\label{eq:1kfn}
	 f\left(y/x \right) \in \mathcal{M}_{pD} \to x^{p} f\left(y /x \right) \in H^{0 } \left(X,L^{p} \right),
	\end{align}
where $ f$ is any polynomial of degree $ \le p$.

	By the second equation of \eqref{eq:tqfx}, for $ p \ge 1 $, the long exact sequence of \eqref{dia:vooe} becomes a short exact sequence.
	Under the identification \eqref{eq:1kfn}, we can rewrite it as 
	\begin{equation}\label{dia:e3gr}
		\begin{tikzcd}
			0 & \mathcal{M}_{\left(p-1\right)D} &	\mathcal{M}_{pD}  & \mathbf{C}  & 0, 
			 \arrow["", from=1-1, to=1-2]
			 \arrow["a", from=1-2, to=1-3]
			 \arrow["b", from=1-3, to=1-4]
			 \arrow["", from=1-4, to=1-5]
		\end{tikzcd} 
	\end{equation}
	where $ a$ is the obvious inclusion, and $ b$ takes the coefficient of the highest order term $ \left(y/x\right)^{p} $.\footnote{A choice of $ s_{D}^{L} $ is required to  identity $ L_{|D} $ with $ \mathbf{C} $.
	We omit the detail.\label{fn:1}
	}
  With the above convention,\footnote{Since $ H^{1} \left(X,\mathbf{C} \right)= 0 $, we take $ s_{0 }^{1} $ to be trivial.} we have  
	\begin{align}\label{eq:wf31}
  s_{p} = 1\wedge \frac{y}{x}   \wedge \ldots \wedge \left(\frac{y}{x} \right)^{p} \in \det \left(\mathcal{M}_{pD}\right). 
	\end{align}

	By Example \ref{exa:akdc}, the associated Slater determinant is given by
	\begin{align}\label{eq:dajo}
	 \prod_{0 \le i<j\le p+1} \left(\frac{y_{j}}{x_{j} }  - \frac{y_{i}}{x_{i} }  \right) \cdot \prod_{i= 0 }^{p} s_{D}^{p}   \left(\left[x_{i} : y_{i} \right]\right) = \prod_{0 \le i<j\le p+1} \left(\frac{y_{j}}{x_{j} }  - \frac{y_{i}}{x_{i} }  \right) \cdot \prod_{i= 0 }^{p} x_{i} ^{p}.
	\end{align}

	Using the embedding $ y \in  \mathbf{C} \to [1:y]\in \mathbf{CP}^{1} $, the above Slater determinant restricts to 
	\begin{align}\label{eq:dajo1}
	 \prod_{0 \le i<j\le p+1} \left(y_{j}  - y_{i}  \right). 
	\end{align}
 Equation \eqref{eq:dajo} is compatible with \cite[p.~84]{Klevtsov16}.

\subsection{The case where $ X$ is an elliptic curve}\label{sXg1}
In this section, we assume that $ X$ is an elliptic curve.
Then, $ g= 1$.

Let $ \mathbf{H} $ be the Poincaré upper half plane.
If $ \tau\in \mathbf{H} $, put 
\begin{align}\label{eq:yesv}
 \Lambda_{\tau } = \mathbf{Z} + \tau \mathbf{Z} \subset \mathbf{C}.
\end{align}
Then, $ \Lambda_{\tau } $ is a lattice in $ \mathbf{C} $.
Up to an action of $ \mathrm{SL}_{2} \left(\mathbf{Z} \right)$ on $ \mathbf{H} $, there exits a unique $ \tau\in \mathbf{H} $ such that 
\begin{align}\label{eq:22yy}
 X= \mathbf{C} / \Lambda_{\tau }.
\end{align}
If $ z\in \mathbf{C} $, we denote by $ [z]$ the corresponding point in $ X$.

Let $ D$ be the reduced divisor associated to $ [0 ]\in X$.
Let $ L$ and $ s_{D} $ be respectively the holomorphic line bundle and the canonical holomorphic section associated to $ D$.\footnote{The objects corresponding to any other reduced divisor of degree one can be obtained through translation.
}
For $ p \ge 1$, we have 
\begin{align}\label{eq:2vyi}
& \dim H^{0 } \left(X,L^{p} \right)= p,&\dim H^{1 } \left(X,L^{p} \right)= 0 .
\end{align}

Let us construct $ L$ and $ s_{D} $ using the Weierstrass functions of $ \wp_{\tau },\zeta_{\tau } , \sigma_{\tau }$.
Following \cite[(III.2.1), (IV.1.1), (IV.2.5) ]{Chandrasekharan85}, for $ z\in \mathbf{C} \setminus \Lambda_{\tau } $,   
\begin{align}\label{eq:bsvq}
 &\wp_{\tau}\left(z\right)= \frac{1}{z^{2} } + \sum_{\omega \in \Lambda_{\tau } \setminus \left\{0\right\} }^{} \left(\frac{1}{\left(z-\omega \right)^{2} } - \frac{1}{\omega^{2} }\right),
 &\zeta_{\tau} \left(z\right)= \frac{1}{z} + \sum_{\omega \in \Lambda_{\tau } \setminus \left\{0\right\} }^{} \left(\frac{1}{z-\omega } + \frac{1}{\omega } + \frac{z}{\omega ^{2} } \right),
\end{align}
and for $ z\in \mathbf{C} $, 
\begin{align}\label{eq:wcm1}
\sigma_{\tau } \left(z\right)&= z \prod_{\omega\in \Lambda_{\tau }\setminus \left\{0\right\}  } \left(1-\frac{z}{\omega } \right)\exp\left(\frac{z}{\omega }  + \frac{1}{2} \frac{z^{2}}{\omega^{2}} \right).
\end{align}

By the first formula of \eqref{eq:bsvq}, $ \wp_{\tau}$ is an even $ \Lambda_{\tau} $-periodic meromorphic function on $ \mathbf{C} $.
Its poles are located at $ \Lambda_{\tau }  $.
As $ z \to 0 $, we have
\begin{align}\label{eq:bkt2}
 \wp_{\tau}\left(z\right) =  z^{-2} + \mathcal{O} \left(1\right).
\end{align}
For $ k \ge 1$, as $ z\to 0 $, the $ k$-th derivation satisfies
\begin{align}\label{eq:hyak}
 \wp_{\tau}^{\left(k\right)}\left(z\right) =  \left(-1\right)^{k} \left(k+1\right)!z^{-2-k}+ \mathcal{O} \left(z^{-k} \right).
\end{align}

Since $ \wp_{\tau}$ has non residue on $ \mathbf{C} $, its primitive exists.
By \eqref{eq:bsvq}, $\zeta_{\tau}$ is the unique odd  primitive of $ -\wp_{\tau}$.
It has simple poles with residue $ 1$ at $ \Lambda_{\tau }$.
Using the fact that $ \wp_{\tau } $ is periodic and that $ \zeta_{\tau } $ is odd, we have 
\begin{align}\label{eq:eixy}
& \zeta_{\tau} \left(z+1\right)= \zeta_{\tau} \left(z\right)+2\zeta_{\tau}\left(1/2\right) ,&\zeta_{\tau} \left(z+\tau \right)= \zeta_{\tau} \left(z\right)+2\zeta_{\tau} \left(\tau/2 \right).
\end{align}
As an application of residue theorem, we get the classical Legendre relation \cite[Theorem IV.2]{Chandrasekharan85}, 
\begin{align}\label{eq:bgwz}
 \zeta_{\tau}\left(1/2\right)\tau -\zeta_{\tau} \left(\tau /2\right)= i \pi. 
\end{align}

Since all the residues of $ \zeta_{\tau} $ are integer $ 1$, the exponential of a primitive of $ \zeta_{\tau}$ exists and is holomorphic with only simple zeros at $ \Lambda_{\tau }$.
By \eqref{eq:wcm1}, $ \sigma_{\tau} $ is such a function with the additional property that it is odd and satisfies, as $ z \to 0 $, 
\begin{align}\label{eq:dizl}
 \sigma_{\tau }  \left(z\right)= z + \mathcal{O} \left(z^{3} \right).
\end{align}
Since $ \sigma_{\tau } $ is odd, by \eqref{eq:eixy}, we have 
\begin{align}\label{eq:fctk}
& \sigma_{\tau }  \left(z+1\right)= - e^{2 \zeta_{\tau} \left(1/2\right) \left(z+ 1/2 \right)}  \sigma_{\tau }  \left(z\right),&\sigma_{\tau }  \left(z+\tau \right)= - e^{2\zeta_{\tau} \left(\tau /2\right) \left(z+ \tau /2 \right)}  \sigma_{\tau } \left(z\right).
\end{align}

The relation \eqref{eq:fctk} indicates that $ \sigma _{\tau } $ is a section of a line bundle on $ X$. 
Indeed, using \eqref{eq:bgwz}, it is easy to see that there is an action of $ \mathbf{Z}^{2} $ on $ \mathbf{C}^{2} $ such that if $ \left(z,\lambda \right)\in \mathbf{C}^{2} $, 
\begin{align}\label{eq:xq2k}
 & \left(1,0 \right)\cdot \left(z,\lambda \right)= \left(z+1,-e^{2 \zeta_{\tau} \left(1/2\right) \left(z+ 1/2 \right)}\lambda \right),&\left(0, 1\right)\cdot \left(z,\lambda \right)= \left(z+\tau, - e^{2\zeta_{\tau} \left(\tau /2\right) \left(z+ \tau /2 \right)}  \lambda \right).
\end{align}
The quotient of $ \mathbf{C}^{2} $ by this $ \mathbf{Z}^{2} $-action defines a holomorphic line bundle on $ X$.
Therefore, $ \sigma_{\tau } $ can be identified with its holomorphic section.
Since $ \sigma_{\tau } $ has only simple zero at $ \Lambda_{\tau} $, we see that the above line bundle and section are just $L$ and $s_{D}$.

Chosen
\begin{align}\label{eq:qvgz}
 s_{1} = s_{D}\in \det H\left(X,L\right),
\end{align}
let us write down a formula for $ s_{p}\in \det H\left(X,L^{p} \right)$ for $ p \ge 2$.\footnote{This is equivalent to taking $ E= L$ and replacing $ p$ by $ p-1$.
The case starting directly from $ s_{0 } $ with $ E$ trivial is left to the reader.
}
 

By \eqref{eq:bkt2} and \eqref{eq:hyak}, for $ p\ge 1$, the space $ \mathcal{M}_{pD} $ is the $ p$-dimensional vector space spanned by the constant function $ 1\in \mathbf{C} $, together with $\wp_{\tau}^{\left(k\right)} $ with $ 0 \le k \le p-2$.
In the current setting, we still have the exact sequence \eqref{dia:e3gr}, where the morphism $ b$ corresponds to taking the coefficient of the term $ \left(-1\right)^{p} \frac{ \wp_{\tau}^{\left(p-2\right)}}{\left(p-1\right)!}$.\footnote{cf. Footnote \ref{fn:1}.}
As in \eqref{eq:wf31},  the above implies
\begin{align}\label{eq:sdwm}
 s_{p} = 1 \wedge \wp_{\tau} \wedge \left(-\frac{\wp_{\tau}^{\left(1\right)}}{2}\right)  \wedge \cdots \wedge\left(\left(-1\right)^{p} \frac{ \wp_{\tau}^{\left(p-2\right)}}{\left(p-1\right)!} \right)\in \det \mathcal{M}_{pD},
\end{align}
Similar to \eqref{eq:dajo},  the Slater determinant associated to \eqref{eq:sdwm} is given by 
\begin{align}\label{eq:3j2n}
 \det \begin{pmatrix}
 1 & 1 & \cdots & 1\\
 \wp_{\tau}\left(z_{1} \right) & \wp_{\tau}\left(z_{2} \right) & \cdots & \wp_{\tau}\left(z_{p} \right)\\
\vdots  & \vdots & \ddots & \vdots\\
 \left(-1\right)^{p} \frac{ \wp_{\tau}^{\left(p-2\right)}\left(z_{1} \right)}{\left(p-1\right)!}  &  \left(-1\right)^{p} \frac{ \wp_{\tau}^{\left(p-2\right)}\left(z_{2} \right)}{\left(p-1\right)!} & \cdots &  \left(-1\right)^{p} \frac{ \wp_{\tau}^{\left(p-2\right)}\left(z_{p} \right)}{\left(p-1\right)!}
 \end{pmatrix}\cdot \prod_{i= 1}^{p} s_{D}^{p} \left(z_{i} \right). 
\end{align}
By \cite[p.~36]{Fay73}, we can rewritten \eqref{eq:3j2n} as
\begin{multline}\label{eq:gpya}
\frac{\sigma_{\tau}\left(z_{1} + \ldots +z_{p} \right) \prod_{1\le i < j\le p} \sigma_{\tau}\left(z_{i} -z_{j} \right)}{\prod_{i= 1}^{p} \sigma_{\tau}^{p}  \left(z_{i} \right)}\cdot \prod_{i= 1}^{p} s_{D}^{p} \left(z_{i} \right)\\
= \sigma_{\tau}\left(z_{1} + \ldots +z_{p} \right) \prod_{1\le i < j\le p} \sigma_{\tau}\left(z_{i} -z_{j} \right).
\end{multline}

Let us compare \eqref{eq:gpya} to \cite[Theorem 3.3]{BurbanKlevtsov24}.
Indeed, in \cite[Section 2]{BurbanKlevtsov24}, the authors use the another line bundle $ L^{\prime } $ obtained by the $ \mathbf{Z}^{2} $-action satisfying
\begin{align}\label{eq:wszd}
& \left(1,0 \right)\cdot \left(z,\lambda \right)= \left(z+1,\lambda \right),&\left(0, 1\right)\cdot \left(z,\lambda \right)= \left(z+\tau, e^{-2i \pi  z - i \pi \tau }  \lambda \right).
\end{align}

We claim that $ L^{\prime } $ is a line bundle associated to the divisor defined by $ \left[\left(1+\tau\right)/2\right]$, or more precisely, if $ T: z \to z + \left(1+\tau\right)/2$ is an automorphism of $ X$, then $ L^{\prime }  \simeq T^{*} L$.
Indeed, if $ s$ is a section of $ L$, then $ s\left(z+\left(1+\tau\right) /2\right)e^{-\zeta_{\tau } \left(1/2\right) \left(z+\left(1+\tau \right)/2\right)^{2} + i \pi \left(z+\left(1+\tau \right)/2\right)}$ defines a section in $ L^{\prime }$.
This gives our claim.

Observe that the above process consists of two steps.
We multiply firstly the section by $ e^{-\zeta_{\tau } \left(1/2\right) z^{2} + i \pi z }$, which does not change the isomorphism class of $ L$.
Then, we apply a shift $ \left(1+\tau\right) /2$ to the variable $ z$, which changes the isomorphism class of $ L$.

Recall that if $ a,b\in \mathbf{R} $, the Jacobi theta function \cite[p.~10]{Mumford83} is defined by
\begin{align}\label{eq:cb2c}
 \theta_{a,b,\tau } \left(z\right)= \sum_{n\in \mathbf{Z} }^{} e^{i \pi \tau \left(n+a\right)^{2} + 2 i \pi \left(n+a\right)\left(z+b\right)}.
\end{align}
Note that $  \theta_{\frac{1}{2} ,\frac{1}{2},\tau } \left(z\right)$ coincides with $ -\theta\left(z,\tau\right)$, a theta function defined in \cite[(V.1.1)]{Chandrasekharan85}.

\begin{prop}\label{prop:ecxa}
 Under the above identification, the section in \eqref{eq:gpya} becomes 
 \begin{multline}\label{eq:rxhw}
 -{i^{p} \exp\left( \frac{i \pi \tau}{4} p^{2}  \right)}{\left[\theta_{\frac{1}{2}, \frac{1}{2},\tau }^{\prime }  \left(0 \right)\right]^{-1-\frac{\left(p-1\right)p}{2} }}\\ \times \theta_{\frac{p-1}{2}, \frac{p-1}{2},\tau  } \left(z_{1} + \ldots +z_{p}  \right) \prod_{1\le i < j\le p}\theta_{\frac{1}{2}, \frac{1}{2},\tau  } \left(z_{i} -z_{j} \right) . 
\end{multline}
\end{prop}
\begin{proof}
Let us identify the section \eqref{eq:gpya} following the above mentioned two steps.

After the first step, our section \eqref{eq:gpya} becomes
\begin{align}\label{eq:ayfn}
 \sigma\left(z_{1} + \ldots +z_{p} \right) \prod_{1\le i < j\le p} \sigma\left(z_{i} -z_{j} \right) \prod_{j= 1}^{p} e^{-p\zeta_{\tau } \left(1/2\right) z^{2}_{j}  + p i \pi z_{j} }. 
\end{align}
By \cite[Theorem V.2]{Chandrasekharan85}, we have 
\begin{align}\label{eq:rxdd}
\sigma_{\tau } \left(z\right)=  \frac{\theta_{\frac{1}{2}, \frac{1}{2},\tau } \left(z\right) }{\theta_{\frac{1}{2}, \frac{1}{2},\tau }^{\prime }  \left(0 \right)} e^{ \zeta_{\tau } \left(1/2\right) z^{2} } .
\end{align}
By \eqref{eq:rxdd}, we can rewrite \eqref{eq:ayfn} as  
\begin{align}\label{eq:ocqq}
 \left[\theta_{\frac{1}{2}, \frac{1}{2},\tau }^{\prime }  \left(0 \right)\right]^{-1-\frac{\left(p-1\right)p}{2} } \theta_{\frac{1}{2}, \frac{1}{2},\tau  } \left(z_{1} + \ldots +z_{p} \right) e^{  p i \pi \sum_{j= 1}^{p} z_{j} }\prod_{1\le i < j\le p}\theta_{\frac{1}{2}, \frac{1}{2},\tau  } \left(z_{i} -z_{j} \right) . 
\end{align}

Using \eqref{eq:cb2c}, it is easy to verify that 
\begin{align}\label{eq:qkxo}
 \theta_{\frac{1}{2}, \frac{1}{2}, \tau } \left(z+\frac{p\left(1+\tau \right)}{2} \right)e^{ p i \pi z}&= \left(-1\right)^{\frac{p\left(p+1\right)}{2} } e^{-\frac{i \pi \tau }{4} p^{2}  }\theta_{\frac{p+1}{2}, \frac{p+1}{2}, \tau } \left(z\right),\\
  \theta_{\frac{p+1}{2}, \frac{p+1}{2},\tau  } \left(z\right)&= \left(-1\right)^{p-1} \theta_{\frac{p-1}{2}, \frac{p-1}{2}, \tau } \left(z\right).
	\notag
\end{align} 
Applying the shift $ \left(1+\tau\right) /2$ to each variable $ z_{j}| _{j= 1,2, \ldots, p} $ in the expression \eqref{eq:ocqq},  and using \eqref{eq:qkxo}, we get \eqref{eq:rxhw}, and finish the proof of our proposition.
\end{proof}

The first line of \eqref{eq:rxhw} is a constant depending only on $ \tau$ and $ p$.
The second line of \eqref{eq:rxhw} coincides with \cite[(3.8)]{BurbanKlevtsov24}.
Therefore, our section \eqref{eq:gpya} is compatible with \cite[Theorem 3.3]{BurbanKlevtsov24}.

\section{The geometric Zabrodin-Wiegmann conjecture}\label{S:ZW}
Given metric data $ \left(\omega^{X}, h^{L}, h^{E} \right)$, we define the partition function $ Z_{p} \left(s_{0 }, s^{L}_{D}, s^{E}_{D} \right)$ to be the square of the $ L^{2} $-norm of the canonical section $ s_{p}\in \lambda_{p} \left(E\right)$. 
The aim of this section is to establish, as $ p\to +\infty$, an asymptotic expansion for $ \log Z_{p} \left(s_{0 }, s^{L}_{D}, s^{E}_{D} \right)$ when $ h^{L} $ has a positive curvature. 
The coefficients of the expansion are explicitly determined up to terms of order $ \mathcal{O} \left(1\right)$.
If $ \left(\omega^{X},h^{L}  \right)$ satisfies the prequantization condition, we can further determine the constant term of the expansion. 

To show our results, we need the corresponding results for the Quillen norm of $ s_{p}  $ and the analytic torsion of $ L^{p} \otimes E$.
The Quillen norm of $ s_{p}$ can be explicitly evaluated and is a polynomial of degree $2$  in $ p$.
This will be detailed in the next section.
The full expansion of the analytic torsion is indeed established by Bismut-Vasserot \cite{BismutVasserot89} and Finski \cite{Finski18}.

This section is organised as follows. 
In Section \ref{sPrequ}, we introduce our metric data and we recall some constructions related to the positive Hermitian line bundle and prequantization conditions.

In Section \ref{sParfun}, we introduce the partition function $ Z_{p} \left(s_{0 }, s_{D}^{L}, s_{D}^{E} \right)$.

In Section \ref{s:State}, we state the main results of our article.

Finally, in Section \ref{s:proof}, we show our main results, while the explicit evaluation of some Quillen norm is deferred to Section \ref{S:EQN}.

In this section, we use the notation of the previous sections. 
In particular, $ L$ is associated to an effective reduced divisor $ D$ and $ s_{D} $ is the canonical holomorphic section of $ L$.
We assume that $ h^{L} $ has a positive curvature.

\subsection{The Chern form and related constructions}\label{sPrequ}
Let $ \omega^{X} $ be a K\"ahler form on $ X$, and let $ h^{TX} $ be the induced Hermitian metric on $ TX$.
Let $ h^{L}, h^{E} $ be Hermitian metrics on $ L$ and $ E$. 

Let $  \nabla^{L} $ be the Chern connection on $   \left(L,h^{L} \right)$.
Let $ R^{L}$ be the corresponding curvature.
Put 
\begin{align}\label{eq:yagj}
c_{1} \left(L,h^{L} \right)= \frac{i}{2\pi } R^{L}.
\end{align}
Then, $ c_{1} \left(L,h^{L} \right)$ is the Chern-Weil representative of the first Chern class $ c_{1} \left(L\right)$.
Similarly, we define $c_{1} \left(TX, \omega^{X} \right), c_{1} \left(E,h^{E} \right)$ to be the first Chern forms of $ \left(TX, h^{TX} \right)$ and $ \left(E,h^{E} \right)$.
More generally, let $ {\rm Td}\left(TX,\omega^{X} \right) $ and  $  {\rm ch}\left(E,h^{E} \right)$ be respectively the Todd form and the Chern character form of $ \left(TX,h^{TX} \right)$ and $ \left(E,h^{E} \right)$.
Since $ \dim X = 1$, we have 
\begin{align}\label{eq:ej3n}
&	{\rm Td}\left(TX,\omega^{X} \right)= 1+\frac{1}{2} c_{1} \left(TX,\omega^{X} \right),& {\rm ch} \left(E,h^{E} \right)= {\rm rk}\left(E\right)+ c_{1} \left(E,h^{E} \right).
\end{align}

In the sequel, we will assume that $\left(L, h^{L}\right) $ is a positive Hermitian line bundle.
Equivalently\footnote{This is only true on Riemann surfaces.}, there is a smooth real function $ r^{L}: X\to \mathbf{R}   $ such that 
\begin{align}\label{eq:ajyc1}
	\frac{i}{2\pi} R^{L} = \exp\left(r^{L} \right)  \omega^{X}.
\end{align}

The pair $ \left(\omega^{X}, h^{L} \right)$ is said to satisfy the prequantization condition, if $ r^{L}= 0  $, i.e.,  
\begin{align}\label{eq:ajyc}
	\frac{i}{2\pi} R^{L} = \omega^{X}. 
\end{align}
In this case, 
\begin{align}\label{eq:sazu}
 \deg L = \int_{X}^{}\omega^{X}.
\end{align}


\subsection{The partition function $ Z_{p} \left(s_{0 }, s_{D}^{L}, s_{D}^{E} \right)$}\label{sParfun}
The metrics $ h^{L}, h^{E} $ induce respectively obvious metrics $ \left| \cdot \right|_{\lambda_{\rm ev} \left(L_{|D} \right)}, \left| \cdot \right|_{\lambda \left(E_{|D} \right)}$ on $ \lambda_{\rm ev} \left(L_{|D} \right)$, $ \lambda \left(E_{|D} \right)$.
Moreover, for $ p\in \mathbf{N} $,  $ \left(h^{L}, h^{E} \right)$ induces a metric $ h^{L^{p} \otimes E} $ on $ L^{p} \otimes E$.
The corresponding analytic torsion will be denoted by $ \tau_{p} \left(\omega^{X}, h^{L}, h^{E} \right)$.
Then, $ \lambda_{p} \left(E\right)$  is equipped with the Quillen metric $ \left\| \cdot \right\|^{{\rm Q}}_{\lambda_{p} \left(E\right)} $ and $ L^{2} $-metric $ \left| \cdot \right|_{\lambda_{p} \left(E\right)} $. 
These two metrics are related by 
\begin{align}\label{eq:h1mg}
 \left\| \cdot \right\|^{{\rm Q}}_{\lambda_{p} \left(E\right)} =  \left| \cdot \right|_{\lambda_{p} \left(E\right)}\exp\left(\tau_{p}\left(\omega^{X}, h^{L}, h^{E} \right) \right).
\end{align}

\begin{defin}\label{def:fnxw}
	Given non-zero data $ s_{0 } \in \lambda \left(E\right)$, $ s_{D}^{L} \in \lambda_{\rm ev} \left(L_{|D} \right)$, and $ s_{D}^{E} \in \lambda \left(E_{|D} \right)$, we define the partition function\footnote{For simplification, we omit the explicit dependence on $ \left(\omega^{X}, h^{L}, h^{E} \right)$.} by 
	\begin{align}\label{eq:m3dl}
	 Z_{p} \left(s_{0}, s_{D}^{L} ,s_{D}^{E}\right)= \left| s_{p} \right|_{\lambda_{p} \left(E\right)}^{2}. 
	\end{align}
\end{defin}

We equip both sides of \eqref{eq:av1q} the metrics induced by $ \left| \cdot \right|_{\lambda_{\rm ev} \left(L_{|D} \right)}, \left| \cdot \right|_{\lambda \left(E_{|D} \right)}$, $ \left\| \cdot \right\|_{\lambda \left(E\right)}^{{\rm Q}} $, $ \left\| \cdot \right\|_{\lambda_{p}  \left(E\right)}^{{\rm Q}} $.
Denote by $ \left\| \sigma_{p}^{0}  \right\|^{{\rm Q}} $ the corresponding norm of $ \sigma_{p}^{0} $.

\begin{prop}\label{prop:suof}
	For $ p\in \mathbf{N} $, the following identity holds, 
	\begin{multline}\label{eq:xqsw}
	 \log Z_{p} \left(s_{0}, s_{D}^{L}  ,s_{D}^{E}  \right)= \log \left\|\sigma_{p}^{0 }  \right\|^{{\rm Q},2}	-2\tau_{p} \left(\omega^{X}, h^{L}, h^{E} \right)+\log \left| s_{0 }  \right|^{2}_{\lambda \left(E\right)}   +2\tau\left(\omega^{X}, h^{E} \right)\\
	 +\frac{p^{2} }{2} {\rm rk}\left(E\right)\log \left|s_{D}^{L}   \right|^{2}_{\lambda_{\rm ev} \left(L_{|D} \right)} +p\left(\frac{1}{2} {\rm rk}\left(E\right)\log \left|s_{D}^{L}   \right|^{2}_{\lambda_{\rm ev} \left(L_{|D} \right)}  +\log \left|s_{D}^{E}   \right|^{2}_{\lambda \left(E_{|D} \right)} \right).
 \end{multline}
 \end{prop}
\begin{proof} 
Our proposition is an immediate consequence of \eqref{eq:bihv}, \eqref{eq:h1mg}, and \eqref{eq:m3dl}.
\end{proof}

\subsection{Statement of our main results}\label{s:State}
Recall that $ s_{D}  $ has simple zeros on $ D$. 
For any connection $ \nabla^{} $ on $ L$, the value $ \nabla^{} s_{D} $ on $ D$ is independent of $ \nabla^{} $, and is a nowhere vanishing section of $ \left(T^{*} X \otimes L\right)_{|D} $ over $ D$.
Denote by $ \partial  s_{D} $ this section on $ D$.
Then, we have an isomorphism of vector bundles on $ D$,  
\begin{align}\label{eq:uvg3}
	\partial s_{D} :  T^{} X_{|D} \to L_{|D}. 
\end{align}

Note that $ \log \left| s_{D}  \right|^{2} $ defines a locally integrable current on $ X$.
Denote by $ \zeta \left(s\right)$ the Riemann zeta function on $ s\in \mathbf{C} $. 
Also, $ \chi\left(X\right)= 2-2g$ is the Euler characteristic of $ X$.

\begin{thm}\label{thm:wlcn}
There exist constants $ \{a_{i} \}_{i\in  2- \mathbf{N} }$ and $ \{b_{i} \}_{i\in 1- \mathbf{N} } $, such that for $ N\in \mathbf{N} $, as $ p\to +\infty$,
 \begin{align}\label{eq:dgid}
	\log Z_{p} \left(s_{0}, s_{D}^{L}, s_{D}^{E} \right)= a_{2} p^{2}  + \sum_{i= 0 }^{N} \left(b_{1-i} \log p +a_{1-i} \right)p^{1-i} +\cO\left(\frac{\log p}{p^{N} }  \right).
 \end{align}
Moreover, 
\begin{align}\label{eq:koff}
	a_{2}  =&  \frac{1}{2}{\rm rk}\left(E\right)\left(\log \left|s_{D}^{L}   \right|^{2}_{\lambda_{\rm ev} \left(L_{|D} \right)} + \int_{X}^{}c_{1} \left(L,h^{L} \right)\log \left| s_{D}  \right|^{2}\right),\\
	b_{1} =&  {  -} \frac{1}{2}{\rm rk}\left(E\right)\deg\left(L\right), \notag\\
	a_{1} =&   \log \left|s_{D}^{E}   \right|^{2}_{\lambda \left(E_{|D} \right)}+\frac{1}{2}{\rm rk}\left(E\right)\Bigg\{\log \left|s_{D}^{L}   \right|^{2}_{\lambda_{\rm ev} \left(L_{|D} \right)}  -\sum_{z \in D} \log \left| \partial s_{D}\left(z\right) \right|^{2} {  -} \int_{X}^{}r^{L}  c_{1} \left(L,h^{L} \right) \Bigg\}\notag\\
	&+\int_{X}^{}{\rm Td}\left(TX,\omega^{X} \right){\rm ch}\left(E,h^{E} \right)\log\left| s_{D} \right|^{2},\notag\\
	b_{0 } =&  {  -} \frac{1}{3}{\rm rk}\left(E\right)\chi\left(X\right) {  -} \frac{1}{2}\deg\left(E\right).\notag
\end{align}
Also, if the prequantization condition \eqref{eq:ajyc} holds, the term $ \int_{X}^{}r^{L} c_{1} \left(L,h^{L} \right) $ in $ a_{1} $ vanishes, and   
 \begin{align}\label{eq:dbkb}
	a_{0 } = \log \left|s_{0 } \right|_{\lambda \left(E\right)} ^{2}	+2\tau\left(\omega^{X}, h^{E} \right)  {  -}  \left(\zeta'\left(-1\right)+\frac{\log\left(2\pi \right)}{12}+\frac{7}{24}\right){\rm rk}\left(E\right)\chi\left(X\right)
	{  -} \frac{1}{2}\deg\left(E\right).
 \end{align}
\end{thm}
\begin{proof}
 The proof of our theorem will be given in Section \ref{s:proof}.
\end{proof}

\begin{re}\label{re:1e1u}
	Note that the constants $a_{i} $ with $ i \neq 0 $, as well as all the $  b_{i}  $ are locally calculable. Similar holds for the constant $ a_{0 } - \log \left| s_{0 }  \right|_{\lambda \left(E\right)}^{2} -2 \tau\left(\omega^{X}, h^{E} \right)$.
	Note also that all the $ b_{i} $ are independent of the metrics.
	See Section \ref{s:proof} for a proof of these facts.
\end{re}


Recall that the torus $ \frac{H^{1} _{{\rm dR}} \left(X,\mathbf{R} \right)}{H^{1} _{{\rm dR}} \left(X,\mathbf{Z} \right)}$ equips a metric induced from the $ L^{2} $-metric and Hodge theory.
Denote by $ {\rm vol}_{L^{2} } \left(\frac{H^{1} _{{\rm dR}} \left(X,\mathbf{R} \right)}{H^{1} _{{\rm dR}} \left(X,\mathbf{Z} \right)}\right)$ the corresponding volume.

\begin{cor}\label{cor:hzng}
	Assume that $ E = \mathbf{C} $ is trivial, $ s_{0 } $ is defined in \eqref{eq:agmh}, and the prequantization condition \eqref{eq:ajyc} holds.
	Then, 
	\begin{align}\label{eq:acy2}
	 a_{2} &= \frac{1}{2}\left(\log \left|s_{D}^{L}   \right|^{2}_{\lambda_{\rm ev} \left(L_{|D} \right)}  +  \int_{X}^{}c_{1} \left(L,h^{L} \right)\log \left| s_{D}  \right|^{2}\right),\\
	 b_{1} &= {  -}  \frac{1}{2}\deg\left(L\right), \notag\\
	 a_{1} &= \frac{1}{2} \left(\log \left|s_{D}^{L}   \right|^{2}_{\lambda_{\rm ev} \left(L_{|D} \right)}  -\sum_{z \in D} \log \left| \partial s_{D}\left(z\right) \right|^{2}+\int_{X}^{}c_{1} \left(TX,\omega^{X} \right)\log\left| s_{D} \right|^{2}\right),\notag \\
	 b_{0 } &=  {  -} \frac{1}{3}\chi\left(X\right), \notag\\
	 a_{0 } &= \log \left\{\frac{\deg L}{2\pi } {\rm vol }^{-1} _{L^{2} } \left(\frac{H^{1} _{{\rm dR}} \left(X,\mathbf{R} \right)}{H^{1} _{{\rm dR}} \left(X,\mathbf{Z} \right)}\right)\right\}	+2\tau\left(\omega^{X}\right)  {  -}  \left(\zeta'\left(-1\right)+\frac{\log\left(2\pi \right)}{12}+\frac{7}{24}\right)\chi\left(X\right).\notag
	\end{align}
\end{cor}
\begin{proof}
By Theorem \ref{thm:wlcn}, we need only to show 
\begin{align}\label{eq:nxun}
 \left| s_{0 }  \right|_{\det H\left(X,\mathbf{C} \right)}^{2}  = \frac{\deg L}{2\pi } {\rm vol }^{-1} _{L^{2} } \left(\frac{H^{1} _{{\rm dR}} \left(X,\mathbf{R} \right)}{H^{1} _{{\rm dR}} \left(X,\mathbf{Z} \right)}\right).
\end{align}

Recall that $ s_{0 }^{0}$ and $ s_{0 }^{1} $ are defined in Section \ref{sEtri}. We have the trivial identity, 
\begin{align}\label{eq:nxun1}
	\left| s_{0 }  \right|^{2} _{\det H\left(X,\mathbf{C} \right)} = \left| s_{0 }^{0 }   \right|^{2} _{\det H^{0 } \left(X,\mathbf{C} \right)} \left| s_{0 }^{1}   \right|_{\det H^{1} \left(X,\mathbf{C} \right)}^{-2}. 
 \end{align}
By \eqref{eq:epzu} and \eqref{eq:sazu}, we have   
\begin{align}\label{eq:v3sr}
	\left| s_{0 }^{0 }   \right|^{2}  _{\det H^{0 } \left(X,\mathbf{C} \right)}= \frac{\deg L}{2\pi }.
\end{align}

By Hodge theory, \eqref{eq:nztx} is an isometry with respect to the $ L^{2} $-metric.
By \eqref{eq:dlig}, 
\begin{align}\label{eq:ewd1}
	\left| s_{0 }^{1 }   \right|^{2} _{\det H^{1 } \left(X,\mathbf{C} \right)}= \left| \tau_{\mathbf{Z} }  \right|_{\det H^{1}_{{\rm dR}} \left(X,\mathbf{R} \right)} = {\rm vol }_{L^{2} } \left(\frac{H^{1} _{{\rm dR}} \left(X,\mathbf{R} \right)}{H^{1} _{{\rm dR}} \left(X,\mathbf{Z} \right)}\right).
\end{align}
From \eqref{eq:nxun1}-\eqref{eq:ewd1}, we obtain \eqref{eq:nxun}, and finish the proof of our corollary.
\end{proof}  
   
\subsection{Proof of our main results}\label{s:proof}
\begin{thm}\label{thm:ngzf}
	For $ p\ge 1$, the following identity holds, 
	\begin{multline}\label{eq:pkqt}
		\log \left\|\sigma_{p}^{0 } \right\|^{{\rm Q},2}   = \frac{p^{2} }{2} {\rm rk}\left(E\right)\int_{X} c_{1} \left(L,h^{L} \right)\log \left| s_{D}  \right|^{2} \\
		+ p\Bigg\{ 
    	\int_{X}^{}{\rm Td}\left(TX,\omega^{X} \right){\rm ch}\left(E,h^{E} \right)\log \left| s_{D}  \right|^{2}-\frac{1}{2}{\rm rk}\left(E\right)\sum_{z \in D} \log \left| \partial s_{D}\left(z\right) \right|^{2} 
		\Bigg\}.
	\end{multline}
\end{thm}
\begin{proof}
 The proof of our theorem will be given in Section \ref{s:pBL}.
\end{proof}

\begin{re}\label{re:dx3j}
 Note that $ \log \left\| \sigma_{p}^{0}  \right\|^{{\rm Q},2}  $ contributes only to the terms $ a_{2} $ and $ a_{1} $.
\end{re}

\begin{re}\label{re:dlu3}
 Theorem \ref{thm:ngzf} holds for general $ h^{L} $.
 See Section \ref{S:EQN}.
\end{re}

\begin{thm}\label{thm:zzjm}
	There exist constants $ \{c_{i} \}_{i\in 1-\mathbf{N} }$ and $ \{d_{i} \}_{i\in   1 - \mathbf{N} } $, such that for $ N\in \mathbf{N} $,  as $ p\to +\infty$, we have 
	\begin{align}\label{eq:kapz}
		2\tau_{p} \left(\omega^{X}, h^{L}, h^{E} \right)= \sum_{i= 0 }^{N} \left(c_{1-i}\log p +d_{1-i}\right) p^{1-i} +\cO\left(\frac{\log p}{p^{N} } \right). 
	\end{align}
	All the constants $ c_{i}, d_{i} $ are locally calculable, and $ c_{i} $ are independent of the  metric data $ \left( \omega^{X}, h^{L}, h^{E} \right)$.
	Moreover, 
	\begin{align}\label{eq:if3r}
		&c_{1} = \frac{1}{2}{\rm rk}\left(E\right)\deg\left(L\right),
		&d_{1} = \frac{1}{2}{\rm rk}\left(E\right)\int_{X}^{} r^{L}  c_{1} \left(L,h^{L} \right),\\
		&c_{0 } = \frac{1}{3}{\rm rk}\left(E\right)\chi\left(X\right)+\frac{1}{2}\deg\left(E\right).&&\notag
	\end{align}
  Also, if the prequantization condition \eqref{eq:ajyc} holds, the above $ d_{1} $ vanishes, and 
\begin{align}\label{eq:nely}
 d_{0 } = \left(\zeta'\left(-1\right)+\frac{\log\left(2\pi\right)}{12}+\frac{7}{24} \right){\rm rk}\left(E\right)\chi\left(X\right)+\frac{1}{2}\deg\left(E\right).
\end{align}
 \end{thm}
\begin{proof}
	Our theorem follows from Bismut-Vasserot \cite[Theorem 8]{BismutVasserot89} and Finski \cite[Theorems 1.1 and 1.3]{Finski18}.
    
	Indeed, in \cite[Theorem 8]{BismutVasserot89}, the asymptotic of $ 2\tau_{p} \left(\omega^{X}, h^{L}, h^{E} \right)$  as $ p\to +\infty$ is established up to $ o\left(p \right)$.
	The coefficients $ c_{1}, d_{1} $ are explicitly computed there, providing the first line of \eqref{eq:if3r}. 
	If \eqref{eq:ajyc} holds, we have $ r^{L} =0 $, which implies $ d_{1} = 0 $.
	
	The full asymptotic expansion	together with the descriptions of $ c_{i}, d_{i} $ are given in \cite[Theorem 1.1]{Finski18}. 
	Under the assumption \eqref{eq:ajyc}, $ c_{0 }, d_{0 } $ are explicitly evaluated in \cite[Theorem 1.3]{Finski18}, giving the last equation of \eqref{eq:if3r} and \eqref{eq:nely}.
	Since $ c_{0 } $ is independent of all metrics, the last equation of \eqref{eq:if3r} holds without the assumption \eqref{eq:ajyc}.
	The proof of our theorem is complete.
\end{proof}
  
\begin{proof}[Proof of Theorem \ref{thm:wlcn} and Remark \ref{re:1e1u}]
	By Proposition \ref{prop:suof}, Theorems \ref{thm:ngzf} and \ref{thm:zzjm}, we get our results.
\end{proof}

\section{Evaluation of the Quillen norm of $ \sigma_{1}^{0 } $}\label{S:EQN}

The aim of this section is to evaluate the Quillen norm of $ \sigma_{1}^{0 } $ with respect to $ \left(\omega^{X}, h^{L}, h^{E} \right)$.
This is achieved in two steps.
First, we introduce an adapted metric $ h^{L}_{1} $ on $ L$ such that the triple $ \left(\omega^{X}, h^{L}_{1}, h^{E} \right)$ satisfies Assumption A of Bismut \cite{B90immersion}.
This allows us to use Bismut-Lebeau's embedding formula \cite{BL91} to evaluate the Quillen norm of $ \sigma_{1}^{0 } $ for $ \left(\omega^{X}, h^{L}_{1} , h^{E} \right)$.
Second, we apply the anomaly formula of Bismut-Gillet-Soulé \cite{BGS3} to compute the  Quillen norm of $ \sigma_{1}^{0 } $ with respect to our initial metrics $ \left(\omega^{X}, h^{L}, h^{E} \right)$.

This section is organised as follows.
In Section \ref{s:admet}, we introduce an adapted metric $ h^{L}_{1} $ on $ L$, so that Assumption A is satisfied.
We evaluate the associate Bott-Chern secondary class $ \widetilde{{\rm ch}} \left(L,h^{L}, h^{L}_{1} \right)$, which appears in the anomaly formula of Bismut-Gillet-Soulé.

In Section \ref{sSupBC}, we review the superconnection formalism.

In Section \ref{s:SBC}, we construct the singular Bott-Chern current associated to $ h_{1}^{L} $, a key term in Bismut-Lebeau's embedding formula.

In Section \ref{sQnorm}, we evaluate the Quillen norm of $ \sigma_{1}^{0 } $.

Finally, in Section \ref{s:pBL}, we show Theorem \ref{thm:ngzf}.

We use the notations and assumptions of the previous sections, with the exception that we no longer require $ h^{L} $ to have positive curvature. 
The results in this section remain valid for an arbitrary $ h^{L} $.

\subsection{An adapted metric on $ L$ and Assumption A}\label{s:admet}

Let us specialize the constructions in \cite[Section 1]{B90immersion} to the embedding $ \iota:D\to X $ and the vector bundle $ \left(L \otimes E\right)_{|D} $ on $ D$.
 
Let $ \xi$ be the degree decreasing\footnote{See \cite[Remark 1.4]{B90immersion} for an explanation of this choice of convention.} complex of holomorphic vector bundles on $ X$ concentrated at degree $ 1$ and $ 0 $, 
\begin{equation}\label{dia:qeml}
	\begin{tikzcd}
		0 &  E &	L \otimes E & 0 .
		 \arrow["", from=1-1, to=1-2]
		 \arrow["s_{D} ", from=1-2, to=1-3]
		 \arrow["", from=1-3, to=1-4]
	\end{tikzcd} 
\end{equation}
By \eqref{dia:pgtl}, the complex of sheaves associated to $ \xi $ provides a resolution for $ \iota_{*} \left(L \otimes E\right)_{|D} $.

Let  $ \pi: TX_{|D} \to D$ be the canonical projection.
Let $ y \in C^{\infty} \left(TX_{|D},\pi^{*}  TX_{|D} \right) $  be the tautological section.
We have a complex of holomorphic vector bundles on the total space of $ TX_{|D} $, 
\begin{equation}\label{dia:qeml2}
	\begin{tikzcd}
		0 & \pi^{*}  \left(E_{|D} \right)   &	\pi^{*} \left(L_{|D} \otimes E_{|D} \right)  & 0 .
		 \arrow["", from=1-1, to=1-2]
		 \arrow["\partial_{y} s_{D} ", from=1-2, to=1-3]
		 \arrow["", from=1-3, to=1-4]
	\end{tikzcd} 
\end{equation}

Using the trivial identification $ TX\otimes T^{*} X= \mathbf{C} $ and the isomorphism \eqref{eq:uvg3}, we can rewrite \eqref{dia:qeml2} as
\begin{equation}\label{dia:qeml1}
	\begin{tikzcd}
		0 & \pi^{*} \left(T^{*} X_{|D}\otimes TX_{|D} \otimes  E_{|D} \right)  &	\pi^{*} \left(TX_{|D} \otimes E_ {|D} \right) & 0 .
		 \arrow["", from=1-1, to=1-2]
		 \arrow["i_{y}  ", from=1-2, to=1-3]
		 \arrow["", from=1-3, to=1-4]
	\end{tikzcd} 
\end{equation}
This is just the Koszul complex $K=   \left(\pi^{*}\left(\Lambda \left(T^{*} X_{|D} \right)\otimes TX_{|D} \otimes E_{|D} \right), i_{y} \right)$.

The complexes \eqref{dia:qeml}, \eqref{dia:qeml1} are exactly the complexes in \cite[(1.2), (1.4)]{B90immersion}, and \eqref{dia:qeml2} is the one described in \cite[Theorem 1.2]{B90immersion}.
The identification between \eqref{dia:qeml2} and \eqref{dia:qeml1} is just the last statement of the above reference.

If we equip \eqref{dia:qeml2} and \eqref{dia:qeml1} the induced metrics from $ \left(\omega^{X}, h^{L}, h^{E} \right)$, this identification is not an isometry.
Let $ h^{L}_{1} $ be another Hermitian metric on $ L$ such that 
\begin{align}\label{eq:2gd2}
 \partial s_{D} : \left(TX,h ^{TX} \right)_{|D}  \to \left(L,h^{L}_{1}  \right)_{|D} 
\end{align}
is an isometry.
Since $ D$ is finite, such a metric $ h^{L}_{1} $ exists always. 

Let $ h^{\xi}_{1} $ and $ h^{K}  $ be the Hermitian metrics on $ \xi $ and $K$ induced respectively by $\left( h^{L}_{1},h^{E} \right) $ and $\left(h^{TX}, h^{E} \right)$.
By \eqref{eq:2gd2}, the complex \eqref{dia:qeml2} endowed with the obvious Hermitian metric induced from $ \left(h^{L}_{1}, h^{E} \right)$ is isometric to $ \left(K,h^{K} \right)$.
This metric compatibility condition is precisely Assumption A of Bismut \cite[Definition 1.5]{B90immersion}. 

The obvious notation associated to $ h^{L}_{1} $ will be added a subscript $1$. 
In particular, if $ \nabla^{L \prime\prime} $ is the holomorphic structure of $ L$, $ \nabla^{L}_{1} = \nabla^{L \prime\prime} + \nabla^{L\prime}_{1} $ denotes the Chern connection of $ h^{L}_{1} $, and $ R^{L}_{1} $ is the corresponding curvature.

If $ \alpha = \alpha^{0} +\alpha^{2} \in \Omega^{\rm even} \left(X,\mathbf{C} \right)$, set
\begin{align}\label{eq:k2tu}
 \varphi \alpha = \alpha^{0 } +\frac{\alpha ^{2}}{2i\pi}.
\end{align}
Then,
\begin{align}\label{eq:fy1z}
	&{\rm ch}\left(L,h^{L} \right)= \varphi \exp\left({- R^{L} }\right) ,& {\rm ch}\left(L,h^{L}_{1}  \right)= \varphi \exp\left({- R^{L}_{1}  }\right) .
\end{align}

Let $ \widetilde {{\rm ch}}\left(L,h^{L}, h_{1}^{L} \right) \in \Omega^{{\rm even}}\left(X,\mathbf{R} \right)/d\Omega^{1}\left(X,\mathbf{R} \right)$ be the Bott-Chern secondary class\footnote{If $ X$ is a complex manifold of arbitrary dimension, then the Bott-Chern secondary class lies in $ \bigoplus_{p= 0 }^{\dim X}  \Omega^{p,p}\left(X,\mathbf{R} \right)/\left(\im \partial +\im \overline{\partial}\right) \cap \Omega^{p,p}\left(X,\mathbf{R} \right). $
When $ \dim X= 1$, 
this space reduces to $ \Omega^{{\rm even}} \left(X,\mathbf{R} \right)/ d\Omega^{1} \left(X,\mathbf{R} \right)$. }
of $ L$ associated to $ h^{L} $ and $ h_{1}^{L} $.
Then, 
\begin{align}\label{eq:tfpl}
	{\rm ch}\left(L,h_{1}^{L}\right)-{\rm ch}\left(L,h^{L}\right)= \frac{\overline{\partial} \partial }{2i\pi} \widetilde {{\rm ch}}\left(L,h^{L}, h_{1}^{L} \right).
\end{align}
 
Let $ \phi \in C^{\infty}\left(X,\mathbf{R} \right)$ such that
\begin{align}\label{eq:ngm3}
	h_{1}^{L} = e^{\phi} h^{L}.
\end{align}
 
\begin{prop}\label{prop:fjaw}
 The following identity holds,
 \begin{align}\label{eq:etaj}
	\widetilde {{\rm ch}}\left(L,h^{L}, h_{1}^{L} \right)=  -\phi - \phi c_{1} \left(L,h^{L} \right)+\frac{\phi\overline{\partial} \partial \phi}{4i\pi} \text{ in }  \Omega^{{\rm even}} \left(X,\mathbf{R} \right)/ d\Omega^{1} \left(X,\mathbf{R} \right). 
 \end{align}
\end{prop}
\begin{proof}
	For $ c \in \mathbf{R}$, write $ h_{c}^{L} =e^{c\phi} h^{L} $. 
	Then, $ h_{0 }^{L} = h^{L} $.
	Let $ R^{L}_{c} $ be the curvature of $ \left(L,h_{c}^{L}\right) $. 
	Then, 
\begin{align}\label{eq:ooke}
 R_{c}^{L}  = R^{L} +c \overline{\partial} \partial \phi.
\end{align}
By \cite[Definition 1.26]{BGS1}, we have 
	\begin{align}\label{eq:csou}
		\widetilde {{\rm ch}}\left( L,h^{L}, h^{L}_{1} \right)=  - \varphi\int_{0 }^{1} \phi \exp\left(-R_{c}^{L} \right)dc.
	\end{align}
	By \eqref{eq:ooke} and \eqref{eq:csou}, we get \eqref{eq:etaj} and finish the proof of our proposition.
\end{proof}

Similar, let $ \nabla^{\xi }_{1}  $ be the Chern connection of the Hermitian complex $ \left(\xi,h^{\xi }_{1} \right)$, and let $ {\rm ch}\left(\xi, h^{\xi }_{1} \right)$ be the corresponding Chern character form.
Clearly, 
\begin{align}\label{eq:msay}
	{\rm ch}\left(\xi, h^{\xi }_{1} \right)= {\rm rk}\left(E\right) c_{1} \left(L,h^{L}_{1} \right).
\end{align}

\subsection{Superconnections and the Chern character forms}\label{sSupBC}
We follow the construction of \cite[Section II]{B90immersion}.
Recall that $ s_{D}: E\to L\otimes E $ is the multiplication by $ s_{D} $ on $ E$.
Let $ \overline{s}^{*} _{D1} $ be its adjoint with respect to $ h^{\xi}_{1} $, i.e., 
\begin{align}\label{eq:gt1j}
	\overline{s}^{*} _{D1} = \left\langle s_{D} , \cdot \right\rangle _{h^{L }_{1}}. 
\end{align}
Here, we adopt the convention that $ \left\langle \cdot , \cdot \right\rangle _{h^{L }_{1}} $ is anti-linear on the first component and linear on the second component.

\begin{defin}\label{def:2hpk}
	For $ u\ge 0 $, set
	\begin{align}\label{eq:xzbo}
	 A_{u} =  \sqrt{u}\begin{pmatrix}
	 0  & \overline{s}_{D1}^{*}  \\
	 s_{D}  & 0 
	 \end{pmatrix} + \nabla^{\xi}_{1}.
	\end{align}
\end{defin}
Then, $ A_{u} $ is a superconnection on the complex $\xi$, and $ A_{u}^{2}  $ is an even element in $ \Omega\left(X,\End\left(\xi \right)\right)$.

\begin{prop}\label{prop:bi3t}
For $ u\ge 0 $, the following identity in $ \Omega \left(X,\End\left(\xi\right)\right)$ holds,
	\begin{align}\label{eq:xqdb}
	 A_{u}^{2} = u\left|s_{D} \right|^{2}_{1}  + \sqrt{u}\begin{pmatrix}
	 0  & \overline{ \left( \nabla^{L\prime}_{1}  s_{D}\right)}_{1}^{*}  \\
	 \nabla^{L\prime}_{1}  s_{D}  & 0 
	 \end{pmatrix} +R^{E}+\begin{pmatrix}
	 0  & 0 \\
	 0  & R^{L}_{1} 
	 \end{pmatrix}. 
	\end{align}
\end{prop}
\begin{proof}
	Our proposition is a consequence of \eqref{eq:xzbo}.
\end{proof}

Note that $ \nabla^{L\prime}_{1} s_{D} $ is a section of $ T^{*} X \otimes L$. 
We use the notation $ \left\langle\nabla^{L\prime}_{1} s_{D} ,  \nabla^{L\prime}_{1} s_{D} \right\rangle _{h^{L}_1}$ defined in an obvious way.
If $ z$ is a holomorphic local coordinate on $ X$, then 
\begin{align}\label{eq:fsmk}
 \left\langle\nabla^{L\prime}_{1} s_{D} ,  \nabla^{L\prime}_{1} s_{D} \right\rangle _{h^{L}_1} = d\overline{z}dz \left\langle\nabla^{L}_{1,\frac{\partial }{\partial z} } s_{D} ,  \nabla^{L}_{1,\frac{\partial }{\partial z} }s_{D} \right\rangle _{h^{L}_1} \in \Omega^{2} \left(X, i \mathbf{R} \right).
\end{align}

Let $ N^{H} $ be the number operator on $ \xi $. 
It acts as multiplication by $ 1$ on $ E$ and by $0 $ on $ L \otimes E$.

\begin{prop}\label{prop:jqh1}
For $ u \ge 0 $, the following identities in $ \Omega^{{\rm even}} \left(X,\mathbf{R} \right)$ hold,  
	\begin{align}\label{eq:2i3b}
		\begin{aligned}
		\varphi{\rm Tr_s}\left[\exp\left(-A_{u}^{2} \right)\right]= & {\rm rk}\left(E\right)\left(c_{1} \left(L,h^{L}_{1}  \right)+\frac{u}{2i\pi}\left\langle \nabla^{L\prime}_{1} s_{D} , \nabla^{L\prime}_{1} s_{D} \right\rangle _{h^{L}_{1 }} \right) e^{-u\left| s_{D}  \right|^{2}_{1}  } ,\\
		\varphi{\rm Tr_s}\left[N^{H} \exp\left(-A_{u}^{2} \right)\right]=& -{\rm rk}\left(E\right)\left(1+\frac{c_{1} \left(E,h^{E} \right)}{{\rm rk}\left(E\right)}-\frac{u}{4i\pi}\left\langle \nabla^{L\prime}_{1} s_{D} , \nabla^{L\prime}_{1} s_{D} \right\rangle _{h^{L}_{1} } \right)e^{-u\left| s_{D}  \right|^{2}_{1}  }.
		\end{aligned}
	\end{align}		 
\end{prop}
\begin{proof}
By \cite[Theorem 1.9]{BGS1}, $ \varphi{\rm Tr_s}\left[\exp\left(-A_{u}^{2} \right)\right]$ and $ \varphi{\rm Tr_s}\left[N^{H} \exp\left(-A_{u}^{2} \right)\right]$ are even forms.

We have 
\begin{align}\label{eq:ywkd}
\overline{ \left( \nabla^{L\prime}_{1}  s_{D}\right)}_{1}^{*}  
	\nabla^{L\prime}_{1}  s_{D} &=  - \left\langle \nabla^{L\prime}_{1} s_{D} , \nabla^{L\prime}_{1} s_{D} \right\rangle _{h^{L}_{1} },\\
	\nabla^{L\prime}_{1}  s_{D}\overline{ \left( \nabla^{L\prime}_{1}  s_{D}\right)}_{1}^{*}  &= \left\langle \nabla^{L\prime}_{1} s_{D} , \nabla^{L\prime}_{1} s_{D} \right\rangle _{h^{L}_{1} }\notag.
\end{align}
Here, the sign in the first equation of \eqref{eq:ywkd} is due to the fact that $ s_{D} $ on the left hand side is an odd operator on $ \xi$. The sign disappears in the second equation of \eqref{eq:ywkd} is because of an extra sign from the anti-commutation relations of holomorphic and anti-holomorphic $ 1$-forms.

By \eqref{eq:xqdb} and \eqref{eq:ywkd}, the even degree part of $ \exp\left(-A_{u}^{2} \right) $ is given by
\begin{align}\label{eq:j2tc}
	 \left(1-R^{E}-\begin{pmatrix}
		0  & 0 \\
		0  & R^{L}_{1}
		\end{pmatrix}  + \frac{u}{2} \left\langle \nabla^{L\prime}_{1} s_{D} , \nabla^{L\prime}_{1} s_{D} \right\rangle _{h^{L}_{1} }  \begin{pmatrix}
		-1 & 0 \\
		0  & 1
		\end{pmatrix}
		 \right)e^{-u \left| s_{D}  \right|^{2}_{1} }.
\end{align}
Now \eqref{eq:2i3b} follows from \eqref{eq:j2tc}. The proof of our proposition is complete.
\end{proof}

We have the double transgression formula\footnote{This is a slight modification of \cite[Theorem 1.15]{BGS1} due to that fact the degree on $ \xi$ is decreasing.
}  \cite[Theorem 2.4]{B90immersion} and the convergence properties \cite[Theorems 3.2 and 4.3]{B90immersion}. 
In our case, all of these can be proved directly from \eqref{eq:2i3b}.


\begin{thm}\label{thm:qiti}
For $ u \ge 0$, the form $ \varphi{\rm Tr_s}\left[\exp\left(-A_{u}^{2} \right)\right]$ is closed. 
If $ u> 0 $, the following identity in $ \Omega^{{\rm even}} \left(X,\mathbf{R} \right)$ holds,
	\begin{align}\label{eq:v3wf}
	 \frac{\partial }{\partial u}\varphi{\rm Tr_s}\left[\exp\left(-A_{u}^{2} \right)\right]= \frac{1}{u}\frac{ \overline{\partial}\partial }{2i\pi}\varphi{\rm Tr_s}\left[N^{H} \exp\left(-A_{u}^{2} \right)\right].
	\end{align}	
\end{thm}
Note that when $ u= 0 $, the form $ \varphi{\rm Tr_s}\left[\exp\left(-A_{u}^{2} \right)\right]$ is just the Chern character form of $ \left(\xi, h^{\xi }_{1} \right)$.
In this case, the first equation of \eqref{eq:2i3b} is just \eqref{eq:msay}.

If $ i\in \{ 0, 1,2\}$, if $ \alpha\in \Omega \left(X,\mathbf{R} \right)$, denote $ \alpha^{\left(i \right)}$ the degree $ i $ component of $ \alpha$.
Let $ \left\| \alpha \right\|_{C^{1} } $ be the $ C^{1} $-norm of $ \alpha$.
Let $ \delta_{D} $ be the Dirac current of $ D$.


\begin{thm}\label{cor:tqqz}
There exists $ C> 0 $ such that  if $ \alpha  \in \Omega \left(X,\mathbf{R} \right)$ and if $ u\ge 1$, we have 
\begin{align}\label{eq:dhsh}
\left|\int_{X}^{}	\varphi{\rm Tr_s}\left[\exp\left(-A_{u}^{2} \right)\right] \alpha -{\rm rk}\left(E\right)\sum_{z \in D} \alpha ^{\left(0 \right)} \left(z\right)\right|&\le \frac{C}{\sqrt{u}}\|\alpha \|_{C^{1} },\\
\left|\int_{X}^{}	\varphi{\rm Tr_s}\left[N^{H} \exp\left(-A_{u}^{2} \right)\right] \alpha   - \frac{1}{2}{\rm rk}\left(E\right)\sum_{z \in D} \alpha ^{\left(0 \right)} \left(z\right)\right|&\le \frac{C}{\sqrt{u}}\|\alpha \|_{C^{1} }.\notag
\end{align}
In particular,  as $ u\to +\infty$, we have the convergences of currents,
	\begin{align}\label{eq:fobi}
	 &\varphi{\rm Tr_s}\left[\exp\left(-A_{u}^{2} \right)\right]\to {\rm rk}\left(E\right)\delta_{D},& \varphi {\rm Tr_s}\left[N^{H} \exp\left(-A_{u}^{2} \right)\right]\to \frac{1}{2} {\rm rk}\left(E\right)\delta_{D}. 
	\end{align}	 
\end{thm}

\subsection{The singular Bott-Chern current}\label{s:SBC}

\begin{defin}\label{def:ycdr}
 For $ s\in \mathbf{C} $ and $ 0 < \Re s < 1/2$, set
 \begin{align}\label{eq:hsbb}
	R\left(\xi,h^{\xi}_{1}  \right)\left(s\right) = \frac{1}{\Gamma \left(s\right)} \int_{0 }^{\infty}u^{s-1} \left(\varphi{\rm Tr_s}\left[N^{H} \exp\left(-A_{u}^{2} \right)\right]-\frac{1}{2} {\rm rk}\left(E\right)\delta_{D}\right)du.
 \end{align} 
\end{defin}
By Theorem \ref{cor:tqqz}, $ R\left(\xi,h^{\xi}_{1} \right)\left(s\right) $ is a current valued holomorphic function on $ 0 < \Re s < 1/2$. 
Moreover, it extends to a meromorphic function on $ s\in \mathbf{C} $ such that $ \Re s<1/2$, which is holomorphic at $ s= 0 $.


\begin{defin}\label{def:etgk}
	The Bott-Chern current of Bismut-Gillet-Soulé \cite[Definition 2.4]{BGS4} is defined by 
 \begin{align}\label{eq:dzc2}
	T\left(\xi,h^{\xi}_{1}  \right) = \frac{\partial }{\partial s}_{|s= 0 } R\left(\xi,h^{\xi}_{1} \right).
 \end{align}
\end{defin}

By construction, we have the following identity of currents, 
\begin{align}\label{eq:dgyd}
	{\rm ch}\left( \xi, h^{\xi }_{1}   \right) + \frac{\overline{\partial} \partial }{2i \pi }T\left(\xi,h^{\xi}_{1}  \right)  = {\rm rk}\left(E\right)\delta_{D}.
\end{align}

The singular support of $ T\left(\xi,h^{\xi}_{1}  \right) $ is given by $ D$. 
By \cite[(3.10) in Theorem 3.3]{BGS4}, $ T\left(\xi,h^{\xi}_{1}  \right) $ is in general not locally integrable and has a singularity of type $ \left| s_{D}  \right|^{-2}_{1}$ near $ D$.

If the square root $ L^{1/2} $ is a well-defined holomorphic line bundle, we have 
\begin{align}\label{eq:urdh}
 {\rm ch}\left(L^{1/2}\otimes E,h^{L^{1/2} \otimes E}_{1}  \right)
 = {\rm rk}\left(E\right)+\frac{1}{2} {\rm rk}\left(E\right)c_{1} \left(L,h^{L}_{1} \right)+c_{1} \left(E,h^{E} \right).
\end{align} 
For ease of notation, even in the case when $ L^{1/2} $ is not well-defined, we still use the abstract notation $ {\rm ch}\left(L^{1/2}\otimes E,h^{L^{1/2} \otimes E}_{1}  \right)$ for the right hand side of \eqref{eq:urdh}.
Clearly, 
\begin{align}\label{eq:gnhf}
  {\rm ch}\left(L^{1/2}\otimes E,h^{L^{1/2} \otimes E}_{1}  \right)=  {\rm ch}\left(L^{1/2},h^{L^{1/2} }_{1}  \right) {\rm ch}\left(E,h^{ E}  \right).
\end{align}
 
\begin{prop}\label{prop:kscx}
Up to adding an exact current, the following identity of currents holds,
\begin{align}\label{eq:dkih}
 	T\left(\xi,h^{\xi}_{1}  \right) = {\rm ch}\left(L^{1/2} \otimes E,h^{L^{1/2} \otimes E}_{1} \right)\log \left| s_{D}  \right|^{2}_{1}.
\end{align}
\end{prop} 
\begin{proof}
Note that  $ T\left(\xi, h^{\xi }_{1} \right)$ is a current of even degree.
Let us show \eqref{eq:dkih} according to its degree.

For $ \alpha \in  \Omega^{2} \left(X,\mathbf{R} \right)$, if $ s\in \mathbf{C} $ such that $  \Re(s)<1$,  we have 
\begin{align}\label{eq:vef1}
 \frac{1}{\Gamma\left(s\right)}\int_{0 }^{\infty}\left(\int_{{X} }^{}\alpha e^{-u\left| s_{D}  \right|_{1}^{2}  } \right)u^{s-1} du = \int_{X}^{} \frac{\alpha }{\left| s_{D}  \right|_{1}^{2s} }. 
\end{align}
By \eqref{eq:vef1}, we get 
\begin{align}\label{eq:wpec}
 \frac{\partial }{\partial s}_{|s= 0 }  \frac{1}{\Gamma\left(s\right)}\int_{0 }^{\infty}\left(\int_{{X} }^{}\alpha e^{-u\left| s_{D}  \right|_{1}^{2}  } \right)u^{s-1} du = -\int_{X}^{}\alpha \log \left| s_{D}  \right|^{2}_{1}.
\end{align}
From \eqref{eq:2i3b} and \eqref{eq:wpec}, we deduce 
\begin{align}\label{eq:mobv}
	T^{\left(0 \right)}\left(\xi,h^{\xi}_{1}  \right)   = {\rm rk}\left(E\right)\log \left| s_{D}  \right|^{2}_{1},
\end{align}
which gives \eqref{eq:dkih} at degree $ 0 $.

Thanks to the current description of $H^{2}_{{\rm dR}} \left(X,\mathbf{R} \right)$ \cite[IV.(6.8)]{DemaillyBook}, up to an exact current, $ T^{\left(2\right)} \left(\xi, h^{\xi }_{1} \right)$ is an element in $ H^{2}_{{\rm dR}} \left(X,\mathbf{R} \right)$.
By \eqref{eq:urdh}, the equation \eqref{eq:dkih} at degree $ 2$ is equivalent to 
\begin{align}\label{eq:o3gc}
	\int_{X}^{}T^{\left(2\right)}\left(\xi,h^{\xi}_{1} \right) = \int_{X}^{} \left(\frac{1}{2}{\rm rk}\left(E\right)c_{1} \left(L,h^{L}_{1}  \right) + c_{1} \left(E,h^{E} \right)\right)\log \left| s_{D}  \right|^{2}_{1}.
 \end{align}

By \eqref{eq:wpec},  we see that the contribution of the second term in the second equation of \eqref{eq:2i3b} in $ T\left(\xi,h^{\xi}_{1} \right)$ is given by the second term on the right hand side of  \eqref{eq:o3gc}.
Thus, to establish \eqref{eq:o3gc}, it remains to show 
 \begin{multline}\label{eq:dajv}
 \frac{\partial }{\partial s}_{|s= 0 }  \frac{1}{\Gamma \left(s\right)} \int_{0 }^{\infty}u^{s-1} \left(\frac{u}{4i \pi} \int_{X}^{}\left\langle \nabla^{L\prime}_{1} s_{D} , \nabla^{L\prime}_{1} s_{D} \right\rangle _{h^{L}_{1} }  e^{-u\left| s_{D}  \right|_{1} ^{2} }-\frac{1}{2} |D|\right)du\\
 =  \frac{1}{2}\int_{X}^{}c_{1} \left(L,h^{L}_{1} \right)\log\left| s_{D}  \right|_{1}^{2}.
\end{multline}

By  the first equation of \eqref{eq:2i3b}, Theorem \ref{thm:qiti}, and the first equation of \eqref{eq:fobi}, 
we know that the class of $ \left(c_{1} \left(L,h^{L}_{1}  \right) + \frac{u}{2i\pi}\left\langle \nabla^{L\prime}_{1} s_{D} , \nabla^{L\prime }_{1} s_{D} \right\rangle _{h^{L}_{1} } \right) e^{-u\left| s_{D}  \right|^{2}_{1}  }$ in $ H^{2}_{{\rm dR}} \left(X,\mathbf{R} \right)$ is independent of $ u\in \mathbf{R}_{+} $, and its integral over $ X$ is given by $ |D|$. 
Therefore, 
\begin{align}\label{eq:id2t}
	\frac{u}{4i \pi} \int_{X}^{}\left\langle \nabla^{L\prime}_{1} s_{D} , \nabla^{L\prime}_{1} s_{D} \right\rangle_{h^{L}_{1} }  e^{-u\left| s_{D}  \right|^{2}_{1}  }-\frac{1}{2} |D|= - \frac{1}{2}  \int_{X}^{}c_{1} \left(L,h^{L}_{1}  \right) e^{-u\left| s_{D}  \right|^{2}_{1}  }.
\end{align}
By \eqref{eq:wpec} and \eqref{eq:id2t},  we deduce \eqref{eq:dajv} and finish the proof \eqref{eq:o3gc}.


The proof of our proposition is complete.
\end{proof}

\begin{prop}\label{prop:1yrn}
	Up to adding an exact current, the following identity of currents holds, 
 \begin{align}\label{eq:hsc2}
	{\rm ch}\left(L^{1/2}, h^{L^{1/2} }_{1} \right)\log \left| s_{D}  \right|^{2}_{1} -{\rm ch}\left(L^{1/2}, h^{L^{1/2} } \right)\log \left| s_{D}  \right|^{2}= - \widetilde {{\rm ch}}\left(L,h^{L}, h_{1}^{L} \right)-\frac{\phi}{2} \delta_{D} . 
 \end{align}
\end{prop}
\begin{proof}
 By \eqref{eq:ngm3}, we have 
 \begin{align}\label{eq:toil}
	{\rm ch}\left(L^{1/2}, h^{L^{1/2} }_{1} \right)\log \left| s_{D}  \right|^{2}_{1}= \left({\rm ch}\left(L^{1/2}, h^{L^{1/2} } \right)-\frac{\overline{\partial} \partial \phi }{4i \pi } \right)\left( \log \left| s_{D}  \right|^{2}+\phi\right).
\end{align}
Therefore, 
\begin{multline}\label{eq:3uxr}
	{\rm ch}\left(L^{1/2}, h^{L^{1/2} }_{1} \right)\log \left| s_{D}  \right|^{2}_{1}-{\rm ch}\left(L^{1/2}, h^{L^{1/2} } \right)\log \left| s_{D}  \right|^{2}\\
	= \phi + \frac{\phi}{2}c_{1} \left(L,h^{L} \right) -\frac{\phi\overline{\partial} \partial \phi }{4i \pi }- \frac{\overline{\partial} \partial \phi }{4i \pi } \log \left| s_{D} 	 \right|^{2}.
\end{multline}

Up to adding an exact current, we have the elementary identity,
\begin{align}\label{eq:e3jw}
	\frac{\overline{\partial} \partial \phi }{2i \pi } \log \left| s_{D} 	 \right|^{2} = \phi \frac{ \overline{\partial} \partial }{2i \pi }\left(\log \left| s_{D} 	 \right|^{2} \right).
\end{align}
By Poincar\'e-Lelong formula \cite[V.(13.2)]{DemaillyBook}, we have
\begin{align}\label{eq:zavb}
	\frac{\overline{\partial} \partial }{2i \pi } \left(\log \left| s_{D} 	 \right|^{2}\right)= \delta_{D} -c_{1} \left(L,h^{L} \right).
\end{align}
From \eqref{eq:etaj}, \eqref{eq:3uxr}-\eqref{eq:zavb}, we get \eqref{eq:hsc2} and finish the proof of our proposition.
\end{proof}

\subsection{The Quillen norm of $ \sigma_{1}^{0 }$}\label{sQnorm}
Recall that the metrics $ \left(\omega^{X}, h^{L}_{1}, h^{E} \right)$ satisfy Assumption A.

\begin{thm}\label{thm:qbsp}
	The following identity holds,
	\begin{align}\label{eq:b1i2}
	 \log \left\|\sigma_{1 }^{0 } \right\|^{{\rm Q},2}_{1}  &= \int_{X}^{}{\rm Td} \left(TX,\omega^{X} \right)T\left(\xi, h^{\xi}_{1}  \right) \\
	 &= \int_{X}^{}{\rm Td} \left(TX,\omega^{X} \right){\rm ch}\left(L^{1/2}\otimes E ,h^{L^{1/2} \otimes E}_{1}  \right) \log \left| s_{D}  \right|^{2}_{1}.\notag
	\end{align}
 \end{thm}
\begin{proof}
Up to a factor\footnote{This adjustment arises because, in \cite{BL91}, the authors consider the inverse of the determinant line.} $ -1$, the first identity in \eqref{eq:b1i2} is a special case of the main result of Bismut-Lebeau \cite[equation (0.5) in Theorem 0.1]{BL91}.
Indeed, that equation expresses $ \log \left\| \sigma_{1}^{0}  \right\|_{1}^{{\rm Q},2}  $ as a sum of four terms.
The first term is precisely the right hand side of \eqref{eq:b1i2}.
The second term vanishes since $ D$ is discrete, so that the normal bundle and the tangent bundle coincide.
The third and last terms vanish by a degree comparison argument and by the fact that $ {\rm ch}\left(\xi \right) = {\rm rk}\left(E\right)c_{1} \left(L\right) $ is of degree $ 2$ and the $ {\rm R}$-genus appearing in these terms is also of degree $ 2$.

The second identity of \eqref{eq:b1i2} is a consequence of \eqref{eq:dkih}.
The proof of our theorem is complete.
\end{proof}

Now we can evaluate the Quillen norm of $ \sigma^{0}_{1}$ for a general triple  $\left(\omega^{X}, h^{L}, h^{E} \right)$.
This gives Theorem \ref{thm:ngzf} for $ p= 1$.

\begin{thm}\label{thm:qbsp2}
	The following identity holds,
	\begin{multline}\label{eq:pkqt1}
		\log \left\|\sigma_{1 }^{0 } \right\|^{{\rm Q},2}
		   = \int_{X}^{} {\rm Td}\left(TX,\omega^{X} \right){\rm ch}\left(L^{1/2} \otimes E,h^{L^{1/2} \otimes E} \right)\log \left| s_{D}  \right|^{2}\\
			 - \frac{1}{2}{\rm rk}\left(E\right)\sum_{z \in D} \log \left| \partial s_{D}\left(z\right) \right|^{2}.
	\end{multline}
 \end{thm}
\begin{proof}
By \eqref{eq:ngm3}, we have the obvious identity 
\begin{align}\label{eq:enhh}
 \log \frac{\left\| \cdot \right\|^{2} _{\lambda_{\rm ev} \left(L_{|D} \right),1} }{\left\| \cdot \right\|_{\lambda_{\rm ev} \left(L_{|D} \right)}^{2} \hfill } = \sum_{z \in D}^{} \phi\left(z\right).
\end{align}

By the anomaly formula of Bismut-Gillet-Soul\'e \cite[Theorem 1.23]{BGS3}, we have
\begin{align}\label{eq:hjaa}
	\log \frac{ \left\|\cdot \right\|^{{\rm Q},2}_{\lambda_{1} \left(E\right)} \hfill}{\left\|\cdot \right\|^{{\rm Q},2}_{\lambda_{1} \left(E\right),1} } =  \int_{X}^{}{\rm Td}\left(TX,\omega^{X} \right){\rm ch}\left(E,h^{E} \right)\widetilde {{\rm ch}}\left( L,h^{L}, h^{ L}_{1} \right).
\end{align}
By \eqref{eq:sqh3}, \eqref{eq:enhh}, and \eqref{eq:hjaa}, we have 
 \begin{align}\label{eq:sqov}
	\log \frac{ \left\|\sigma_{1 }^{0 } \right\|^{{\rm Q},2}}{ \left\|\sigma_{1 }^{0 } \right\|^{{\rm Q},2}_{1} } = \int_{X}^{}{\rm Td}\left(TX,\omega^{X} \right){\rm ch}\left(E,h^{E} \right)\widetilde {{\rm ch}}\left( L,h^{L}, h^{ L}_{1} \right)
	+{\rm rk}\left(E\right)\sum_{z \in D}^{} \phi\left(z\right).
 \end{align}
 
By \eqref{eq:hsc2}, we get
\begin{multline}\label{eq:xo2g}
	\int_{X}^{}{\rm Td}\left(TX,\omega^{X} \right){\rm ch}\left(E,h^{E} \right)\widetilde {{\rm ch}}\left( L,h^{L}, h^{ L}_{1} \right)= -\frac{1}{2} {\rm rk}\left(E\right) \sum_{z\in D}^{} \phi\left(z\right)\\
	+ \int_{X}^{}{\rm Td}\left(TX,\omega^{X} \right){\rm ch}\left(E,h^{E} \right)\left\{{\rm ch}\left(L^{1/2} ,h^{L^{1/2} }\right)\log \left| s_{D}  \right|^{2}-{\rm ch}\left(L^{1/2} ,h^{L^{1/2}  }_{1} \right)\log \left| s_{D}  \right|^{2}_{1} \right\}.
\end{multline}

By \eqref{eq:b1i2}, \eqref{eq:sqov}, and \eqref{eq:xo2g}, we get 
\begin{align}\label{eq:snqh}
	\log \left\|\sigma_{1 }^{0 } \right\|^{{\rm Q},2}
	= \int_{X}^{} {\rm Td}\left(TX,\omega^{X} \right){\rm ch}\left(L^{1/2} \otimes E,h^{L^{1/2} \otimes E} \right)\log \left| s_{D}  \right|^{2}+\frac{1}{2} {\rm rk}\left(E\right)\sum_{z\in D}^{} \phi\left(z\right).
\end{align}

Note that by \eqref{eq:2gd2}, if $ z \in  D$, 
\begin{align}\label{eq:jdut}
	 \phi\left(z\right)= - \log \left| \partial s_{D}\left(z\right) \right|^{2}.
\end{align}
By \eqref{eq:snqh} and \eqref{eq:jdut}, we get \eqref{eq:pkqt1}, and finish the proof of our theorem.
\end{proof}

\subsection{Proof of Theorem \ref{thm:ngzf}}\label{s:pBL}
For $ i\ge 1$, replacing $ \left(E,h^{E} \right)$ by $\left(L^{i-1} \otimes E, h^{L^{i-1} \otimes E} \right)$ in  Theorem \ref{thm:qbsp2}, we get 	
\begin{multline}\label{eq:pkqt3}
	\log \left\|\sigma_{i }^{i-1} \right\|^{{\rm Q},2}
		 = \int_{X}^{} {\rm Td}\left(TX,\omega^{X} \right){\rm ch}\left(L^{i-1/2} ,h^{L^{i-1/2} } \right){\rm ch}\left(E,h^{E} \right)\log \left| s_{D}  \right|^{2}\\
		 - \frac{1}{2}{\rm rk}\left(E\right)\sum_{z \in D} \log \left| \partial s_{D}\left(z\right) \right|^{2}.
\end{multline}
Taking a sum of \eqref{eq:pkqt3} over $1\le i \le p$, and using 
\begin{align}\label{eq:bsqp}
	\sum_{i= 1}^{p} {\rm ch}\left(L^{i-1/2} ,h^{L^{i-1/2} } \right)= p+\frac{p^{2} }{2} c_{1} \left(L,h^{L} \right),
\end{align}
we get \eqref{eq:pkqt} and finish the proof of our theorem.  
\qed

  

\bibliographystyle{amsalpha}

\def\cprime{$'$}
\providecommand{\bysame}{\leavevmode\hbox to3em{\hrulefill}\thinspace}
\providecommand{\MR}{\relax\ifhmode\unskip\space\fi MR }
\providecommand{\MRhref}[2]{%
  \href{http://www.ams.org/mathscinet-getitem?mr=#1}{#2}
}
\providecommand{\href}[2]{#2}

\end{document}